\documentclass[12pt]{article}
\usepackage{amsmath, amssymb}
\begin{document}
\centerline{\bf Multiplier Ideal Sheaves in Complex and Algebraic
Geometry}

\bigskip
\centerline{Yum-Tong Siu\ %
\footnote{Partially supported by a grant from the National Science
Foundation.} }

\bigskip  This article is an
expanded version of the talk I gave on August 23, 2004 in the
International Conference on Several Complex Variables in Capital
Normal University, Beijing, China, and will appear in {\it Science
in China}, Series A Mathematics 2005 Volume 48, as part of the
proceedings of the conference. There are two parts in this
article. The first part, which is the main part of the article,
discusses the application, by the method of multiplier ideal
sheaves, of analysis to complex algebraic geometry.  The second
part discusses the other direction which is the application of
complex algebraic geometry to analysis, mainly to problems of
estimates and subellipticity for the $\bar\partial$ operator.

\bigbreak\noindent{\bf Part I. Application of Analysis to
Algebraic Geometry.}

\bigbreak For the application of analysis to algebraic geometry.  We
will start out with the general technique of reducing problems in
algebraic geometry to problems in $\bar\partial$ estimates for Stein
domains spread over ${\mathbb C}^n$.  The $L^2$ estimates of
$\bar\partial$ corresponds to the algebraic notion of multiplier
ideal sheaves. This method of multiplier ideal sheaves has been
successfully applied to effective problems in algebraic geometry
such as problems related to the Fujita conjecture and the effective
Matsusaka big theorem.  It has also been applied to solve the
conjecture on deformational invariance of plurigenera.  There are
indications that it might possibly be used to give a solution of the
conjecture on the finite generation of the canonical ring. Since the
application of the method of multiplier ideals to effective problems
in algebraic geometry are better known, we will only very briefly
mention such applications.  We will explain more the application to
the deformational invariance of plurigenera and discuss the
techniques and ideas whose detailed implementations may lead to a
solution of the conjecture of the finite generation of the canonical
ring.

\bigbreak\noindent{\sc (1.1) Algebraic Geometric Problems Reduced
to $L^{2}$ estimate for domain (spread over $\mathbb{C}^{n}$)}.

\bigbreak The following simple procedure, of removing certain
hypersurfaces and keeping $L^2$ estimates, reduces problems in
algebraic geometry to problems in $\bar\partial$ estimates for Stein
domains spread over ${\mathbb C}^n$.

\medbreak Let $X$ be an $n$-dimensional complex manifold inside
$\mathbb{P}_{N}$.  Let $S$ be some linear ${\mathbb P}_{N-n-1}$
inside ${\mathbb P}_N$ which is disjoint from $X$. We will use $S$
as the light source for a projection.  Let $T$ be some linear
${\mathbb P}_n$ inside ${\mathbb P}_N$ which is disjoint from $S$.
We will use $T$ as the target for a projection.  We define a
projection $\pi:X\to T$ as follows.  For $x\in X$ we define $\pi
(x)\in T$ as the point of intersection of $T$ with the linear
${\mathbb P}_{N-n}$ in ${\mathbb P}_n$ which contains $S$ and $x$.
Then $\pi:X\to T$ makes $X$ a branched cover over $T$.

\medbreak Let $L$ be a holomorphic line bundle over $X$ and $s$ be
a global (non identically zero) meromorphic section of $L$ over
$X$ with pole-set $A$ and zero-set $B$.

\medbreak Let $Z$ be some hypersurface inside $T$ which contains
the infinity hyperplane of $T$ and contains $\pi\left(A\cup
B\right)$ so that $\pi :X-\pi ^{-1}(Z)\rightarrow T-Z$ is a local
biholomorphism.  After we identify ${\mathbb C}^n$ with $T$ minus
the infinity hyperplane of $T$, $X-\pi^{-1}(Z)$ becomes a Stein
domain spread over ${\mathbb C}^n$.

\medbreak Take a metric $e^{-\varphi}$ of $L$ with $\varphi$
locally bounded from above.  For any open subset $\Omega$ of $X$
and any holomorphic function $f$ on $\Omega-\pi^{-1}(Z)$ with
$$\int_{\Omega-\pi^{-1}(Z)}\left|f\right|^2e^{2\log|s|-\varphi}<\infty,$$
the section $fs$ of $L$ can be extended to a holomorphic section
of $L$ over $\Omega$.

\medbreak Algebraic problems concerning $X$ and $L$ involving
sections and cohomology can be translated, through this procedure,
to problems concerning functions and forms on $\Omega$ involving
$L^2$ estimates of $\bar\partial$ for the weight function
$e^{2\log|s|-\varphi}$.

\bigbreak\noindent{\sc (1.2) Multiplier Ideal Sheaves and Effective
Problems in Algebraic Geometry.}

\bigbreak One important concept to facilitate the translation
between algebraic geometry and analysis is that of {\it multiplier
ideal sheaves}. For a plurisubharmonic function $\varphi$ on an open
subset $U$ of ${\mathbb C}^n$ the multiplier ideal sheaf ${\mathcal
I}_\varphi$ is defined as the sheaf of germs of holomorphic
function-germs $f$ on $U$ such that $\left|f\right|^2e^{-\varphi}$
is locally integrable [Nad89, De93]. This concept helps to translate
$L^2$ estimates into algebraic conditions.  For a holomorphic line
bundle $L$ over a compact complex manifold $X$ with a (possibly
singular) metric $e^{-\psi}$ defined by a local plurisubharmonic
function $\psi$, the multiplier ideal sheaf ${\mathcal I}_\psi$ is a
coherent ideal sheaf on $X$.

\medbreak Later in Part II we will discuss another kind of
multiplier ideals [Ko79] and modules which arise from formulating,
in terms of algebraic conditions, the problems of subelliptic
estimates for smooth weakly pseudoconvex domains.

\medbreak Multiplier ideal sheaves ${\mathcal I}_\psi$ have been
used to successfully solve, or make good progress toward the
solution of, a number of algebraic geometric problems such as the
Fujita conjecture, the effective Matsusaka big theorem, and the
deformational invariance of the plurigenera.  Since the use of
multiplier ideal sheaves in effective problems in algebraic geometry
has a somewhat longer history and is better known, our discussion of
the Fujita conjecture and the effective Matsusaka big theorem will
be very brief. We will focus on the deformational invariance of the
plurigenera and the problem of finite generation of the canonical
ring which is related to it and in a certain sense motivates it.

\medbreak\noindent(1.2.1) {\it Fujita Conjecture.} The Fujita
conjecture [Fu87] states that, if $X$ is a compact complex
algebraic manifold of complex dimension $n$ and $L$ is an ample
holomorphic line bundle on $X$, then $mL+K_{X}$ is globally free
for all $m\geq n+1$ and $mL+K_{X}$ is very ample for all $m\geq
n+2$.

\medbreak For the first part of global freeness, the conjecture is
proved for $n=2$ by Reider [Re88], $n=3$ by Ein-Lazarsfeld [EL93,
Fu93], $n=4$ by Kawamata [Ka97], general $n$ with weaker $m\geq
(1+2+\cdots +n)+1$ (which is of order $n^2$) instead of $m\geq n+1$
by Angehrn-Siu [AS95], improved by Helmke [Hel97, Hel99] and by
Heier [Hei02] to $m$ of order $n^{\frac{4}{3}}$.

\medbreak By using the method of higher-order multiplier ideal
sheaves the very ampleness part can be proved for $m\geq m_n$ with
some explicit effective $m_n$ depending on $n$. The detailed
argument for $n=2$ was given in [Si01] which can be modified for the
case of general dimension $n$.  For the much simpler problem of the
very ampleness of $mL+2K_X$, a bound for $m$ of the order $3^n$ is
easily obtained by using the multiplier ideal sheaves ${\mathcal
I}_\psi$ and the fundamental theorem of algebra [De93, De96a, De96b,
ELN94, Si94, Si96a, Si96b].

\medbreak\noindent(1.2.2) {\it Effective Matsusaka Big Theorem.} For
the effective Matsusaka big theorem for general dimension, the best
result up to this point is the following [Si02b] (for earlier
results see [Si93, De96a] and for a more precise bound in dimension
$2$ see [FdB96]).

\medbreak Let $X$ be a compact complex manifold of complex dimension
$n$ and $L$ be an ample line bundle over $X$ and $B$ be a
numerically effective line bundle over $X$. Then $mL-B$ is very
ample for $m$ no less than
$$
C_n\left(L^{n-1}\tilde K_X \right)^{2^{\max(n-2,0)}}
\left(1+\frac{L^{n-1}\tilde K_X}{L^n} \right)^{2^{\max(n-2,0)}},
$$
where
$$
C_n=2^{n-1+2^{n-1}}\left(\prod_{k=1}^n\left(k
2^{\frac{(n-k-1)(n-k)}{2}}\right)^{2^{\max(k-2,0)}}\right)
$$ and $\tilde K_X=\left(2n
{3n-1\choose n}+2n+1\right)L+B+2K_X$.

\bigbreak\noindent{\sc (1.3) Background of Finite Generation of
Canonical Ring.} We now very briefly present the background for
the problem of the finite generation of the canonical ring and its
relation to the deformational invariance of plurigenera.

\bigbreak\noindent(1.3.1) {\it Pluricanonical Bundle.} Let $X$ be
a complex manifold of complex dimension $n$. Let $K_X$ be the
canonical line bundle so that local holomorphic sections of $K_X$
are local holomorphic $n$-forms.   A local holomorphic section $s$
of $K_X^{\otimes m}$ over $X$ is locally of the form
$f\left(dz_1\wedge\cdots\wedge dz_n\right)^m$, where $f$ is a
local holomorphic function and $z_1,\cdots,z_n$ are local
coordinates of $X$.   We will use the additive notation $mK_X$ for
$K_X^{\otimes m}$.

\medbreak\noindent(1.3.2) {\it Blowup of a Point.}  We can blow up
the origin $0$ of ${\mathbb C}^n$ to form $\widetilde{{\mathbb
C}^n}$ which is the topological closure of the graph of the map
${\mathbb C}^n-\left\{0\right\}\to{\mathbb P}_{n-1}$ defined by
$\left(z_1,\cdots,z_n\right)\mapsto\left[z_1,\cdots,z_n\right]$. The
projection $\widetilde{{\mathbb C}^n}\to{\mathbb C}^n$ is from the
natural projection of the graph onto the domain.

\medbreak\noindent(1.3.3) {\it Monoidal Transformation.}  We can do
the blow-up with a parameter space ${\mathbb C}^k$ by using the
product ${\mathbb C}^n\times{\mathbb C}^k$ and blowing up
$\left\{0\right\}\times{\mathbb C}^k$.

\medbreak For a manifold $X$ and a submanifold $D$ we can blow up
$D$ to get another manifold $\tilde X$, because locally the pair
$\left(X,D\right)$ is the same as the pair $\left({\mathbb
C}^n\times{\mathbb C}^k,\,\left\{0\right\}\times{\mathbb
C}^k\right)$ (if $\dim_{\mathbb C}X=n+k$ and $\dim_{\mathbb C}D=k$).
This blow-up is called the {\it monoidal transformation} of $X$ with
nonsingular center $D$.

\bigbreak\noindent(1.3.4) {\it Resolution of Singularities.}
Hironaka [Hi64] resolved singularity of a subvariety $V$ of a
compact complex manifold $X$ by a finite number of successive
monoidal transformations with nonsingular center so that the
pullback of $V$ to the final blowup manifold $\tilde X$ becomes a
finite number of nonsingular hypersurfaces in normal crossing ({\it
i.e.}, they are locally like a subcollection of coordinate
hyperplanes).

\medbreak\noindent(1.3.5) {\it Space of Pluricanonical Sections
Unchanged in Blowup and Blowdown.} An important property of the
pluricanonical line bundle $mK_X$ of a compact complex manifold
$X$ is that global holomorphic pluricanonical sections ({\it
i.e.,} elements of $\Gamma\left(X,mK_X\right)$) remain global
holomorphic pluricanonical sections in the process of blowing up
and blowing down.

\medbreak  More generally, if we have a holomorphic map $\pi:Y\to X$
between two compact complex manifolds of the same dimension and a
subvariety $Z$ of codimension $\geq 2$ in $X$ such that $f$ maps
$Y-f^{-1}(Z)$ biholomorphically onto $X-Z$, then every element $s$
of $\Gamma\left(Y, mK_Y\right)$ comes from the pullback of some
element $s^\prime$ of $\Gamma\left(X,mK_X\right)$, because the
pushforward of $s|_{Y-f^{-1}(Z)}$ is a holomorphic section of the
holomorphic line bundle $mK_X$ over $X-Z$ and can be extended across
the subvariety $Z$ of codimension $\geq 2$ to give a holomorphic
section $s^\prime$ of $mK_X$ over all of $X$.

\medbreak\noindent(1.3.6) {\it Canonical Ring.}  The ring (known
as the canonical ring)
$$
R\left(X,K_X\right)=\bigoplus_{m=0}^\infty\Gamma\left(X,mK_X\right)
$$
is invariant under blow-ups and blow-downs.

\medbreak Two compact projective algebraic complex manifolds related
by blowing-ups and blowing-downs clearly have the same field of
meromorphic functions ({\it i.e.,} are birationally equivalent).

\medbreak In order to get a representative in a birationally
equivalence class which is easier to study, a most important
question in algebraic geometry is the existence of good
representatives called {\it minimal models}.

\medbreak For simplicity, let us focus on complex manifolds which
are of general type.  A complex manifold $X$ of complex dimension
$n$ is {\it of general type} if $\dim_{\mathbb
C}\Gamma\left(X,mK_X\right)\geq cm^n$ for some $c>0$ and for all $m$
sufficiently large.

\bigbreak\noindent(1.3.7) {\it Conjecture on Finite Generation of
the Canonical Ring}. Let $X$ be a compact complex manifold of
general type.  Let $$ R(X,K_{X})=\bigoplus\limits_{m=0}^{\infty
}\Gamma (X,mK_{X})\;.
$$
Then the ring $R(X,K_X)$ is finitely generated.

\bigbreak If the canonical ring $R\left(X, K_X\right)$ is finitely
generated by elements $s_1,\cdots,s_N$ with
$s_j\in\Gamma\left(X,m_jK_X\right)$.  Let $m_0=\max_{1\leq j\leq
N}m_j$. We can use a basis of $\Gamma\left(X,(m_0!)K_X\right)$ to
define a rational map.

\medbreak The image $Y$ may not be regular.  We expect its canonical
line bundle to behave somewhat like that of a manifold.  More
precisely, its canonical line bundle $K_Y$ is expected to satisfy
the following two conditions.

\begin{itemize}
\item[(i)] $K_Y$ (defined from the extension of the canonical line bundle of
the regular part of $Y$) is a ${\mathbb Q}$-Cartier divisor ({\it
i.e.,} some positive integral power is a line bundle) and is
numerically effective.

\item[(ii)] There exists a resolution of singularity $\pi: \tilde
Y\to Y$ such that $K_{\tilde Y}=\pi^*K_Y+\sum_j a_j E_j$, where
$\tilde Y$ is regular and $\left\{E_j\right\}_j$ is a collection of
hypersurfaces of $Y$ in normal crossing, and $a_j>0$ and
$a_j\in{\mathbb Q}$.\end{itemize} A compact complex variety $Y$ of
general type satisfying (i) and (ii) is called a {\it minimal model}
(see [Ka85]).  (For the purpose of comparing (ii) with the manifold
case, we note that, for any proper surjective holomorphic map
$\sigma:\tilde Z\to Z$ of complex manifolds of the same dimension,
$K_{\tilde Z}=\sigma^*K_Z+\sum_j b_j F_j$, where $b_j$ is a positive
integer and each $F_j$ is a hypersurface of $Y$\,.)

\bigbreak\noindent(1.3.8) {\it Minimal Model Conjecture (for
General Type)} Any compact complex manifold $X$ of general type is
birational to a minimal model.

\bigbreak The conjecture is known for general threefolds [Mo88] and
is still open for general dimension. Analysis offers the possibility
of new tools to handle the conjecture.

\medbreak Some consequences of the conjecture have already been
handled with success by new tools in analysis.  A prominent example
is the proof of the deformational invariance of plurigenera for the
algebraic case.  We will first look at the deformational invariance
of plurigenera and its relation to the minimal model conjecture.
Then we will return to the problem of the finite generation of the
canonical ring later.

\bigbreak\noindent{\sc (1.4) Deformational Invariance of
Plurigenera and the Two Ingredients for its Proof.}

\bigbreak The most general form of the conjecture on the
deformational invariance of plurigenera is for the K\"ahler case
which is stated as follows.

\bigbreak\noindent(1.4.1) {\it Conjecture on Deformational
Invariance of Plurigenera.} Let $\pi:X\to\Delta$ be a holomorphic
family of compact complex K\"haler manifolds over the unit $1$-disk
$\Delta\subset{\mathbb C}$.  Let $X_t=\pi^{-1}(t)$ for $t\in\Delta$.
Then $\dim_{\mathbb C}\Gamma\left(X_t, mK_{X_t}\right)$ is
independent of $t$ for any positive integer $m$.

\bigbreak\noindent By the semi-continuity of $\dim_{\mathbb
C}\Gamma\left(X_t, mK_{X_t}\right)$ as a function of $t$, the main
problem is the extension of an element of $\Gamma\left(X_0,
mK_{X_0}\right)$ to an element of $\Gamma\left(X, mK_X\right)$.
Consider the short exact sequence
$$
0\to{\mathcal
O}_X\left(mK_X\right)\stackrel{\Theta}{\longrightarrow}{\mathcal
O}_X\left(mK_X\right)\to{\mathcal O}_{X_0}\left(m K_{X_0}\right)\to
0,$$ where $\Theta$ is defined by multiplication by $t$.  From its
long exact cohomology sequence the vanishing of
$H^1\left(X,{\mathcal O}_X\left(m K_X\right)\right)$ would give the
extension.

\medbreak Let us first look at the case of general type. If we have
a parametrized version of the minimal model conjecture, then we have
the answer [Nak96], because pluricanonical sections are independent
of blowing-ups and blowing-downs and the minimal models have
numerically effective canonical line bundles for which the theorem
of Kawamata-Viehweg [Ka82, Vi82] would still hold with the kind of
mild singularities of the minimal models.

\medbreak While the minimal model conjecture and the K\"ahler case
of the deformational invariance of plurigenera are both still open
(see [Le83, Le85] for some partial results in the K\"ahler case),
the algebraic case has been proved [Si98, Si02a] by using a
``two-tower argument'' and the following two ingredients [Si98,
p.664, Prop.1 and p.666, Prop.2].

\bigbreak\noindent(1.4.2) {\it Global Generation of Multiplier Ideal
Sheaves (Ingredient One)}. Let $L$ be a holomorphic line bundle over
an $n$-dimensional compact complex manifold $Y$ with a metric which
is locally of the form $e^{-\xi}$ with $\xi$ plurisubharmonic. Let
${\cal I}_\xi$ be the multiplier ideal sheaf of the metric
$e^{-\xi}$. Let $A$ be an ample holomorphic line bundle over $Y$
such that for every point $P$ of $Y$ there are a finite number of
elements of $\Gamma(Y,A)$ which all vanish to order at least $n+1$
at $P$ and which do not simultaneously vanish outside $P$. Then
$\Gamma(Y,{\cal I}_\xi\otimes(L+A+K_Y))$ generates ${\cal
I}_\xi\otimes(L+A+K_Y)$ at every point of $Y$.

\bigbreak\noindent(1.4.3) {\it Extension Theorem of Ohsawa-Takegoshi
Type (Ingredient Two).} Let $Y$ be a complex manifold of complex
dimension $n$. Let $w$ be a bounded holomorphic function on $Y$ with
nonsingular zero-set $Z$ so that $dw$ is nonzero at every point of
$Z$. Let $L$ be a holomorphic line bundle over $Y$ with a (possibly
singular) metric $e^{-\kappa}$ whose curvature current is
semipositive. Assume that $Y$ is projective algebraic (or, more
generally, assume that there exists a hypersurface $V$ in $Y$ such
that $V\cap Z$ is a subvariety of codimension at least $1$ in $Z$
and $Y-V$ is the union of a sequence of Stein subdomains
$\Omega_\nu$ of smooth boundary and $\Omega_\nu$ is relatively
compact in $\Omega_{\nu+1}$). If $f$ is an $L$-valued holomorphic
$(n-1)$-form on $Z$ with $$
\int_Z\left|f\right|^2\,e^{-\kappa}<\infty\; ,
$$
then $fdw$ can be extended to an $L$-valued holomorphic $n$-form
$F$ on $Y$ such that
$$
\int_Y\left|F\right|^2\,e^{-\kappa}\leq 8\pi\,e\,\sqrt{2+{1\over
e}}\,\left(\sup_Y\left|w\right|^2\right)
\int_Z\left|f\right|^2\,e^{-\kappa}\; .
$$

\bigbreak\noindent For the first ingredient we actually need its
effective version which is as follows.

\bigbreak\noindent(1.4.4) {\it Theorem (Effective Version of Global
Generation of Multiplier Ideal Sheaves)}. Let $L$ be a holomorphic
line bundle over an $n$-dimensional compact complex manifold $Y$
with a metric which is locally of the form $e^{-\xi}$ with $\xi$
plurisubharmonic.  Assume that for every point $P_0$ of $Y$ one has
a coordinate chart $\tilde U_{P_0}=\{|z^{(P_0)}|<2\}$ of $Y$ with
coordinates
$$
z^{(P_0)}=\left(z^{(P_0)}_1,\cdots,z^{(P_0)}_n\right)
$$
centered at $P_0$ such that the set $U_{P_0}$ of points of $\tilde
U_{P_0}$ where $\left|z^{(P_0)}\right|<1$ is relatively compact in
$\tilde U_{P_0}$. Let $\omega_0$ be a K\"ahler form of $Y$. Let
$C_Y$ be a positive number such that the supremum norm of $d
z^{(P_0)}_j$ with respect to $\omega_0$ is $\leq C_Y$ on $U_{P_0}$
for $1\leq j\leq n$. Let $0<r_1<r_2\leq 1$. Let $A$ be an ample
line bundle over $Y$ with a smooth metric $h_A$ of positive
curvature. Assume that, for every point $P_0$ of $Y$, there exists
a singular metric $h_{A,P_0}$ of $A$, whose curvature current
dominates $c_A\omega_0$ for some positive constant $c_A$, such
that
$$
{h_A\over|z^{(P_0)}|^{2(n+1)}}\leq h_{A,P_0}
$$
on $U_{P_0}$ and
$$
\sup_{r_1\leq|z^{(P_0)}|\leq r_2}{h_{A,P_0}(z^{(P_0)})\over
h_A(z^{(P_0)})}\leq C_{r_1,r_2}
$$
and
$$
\sup_Y{h_A\over h_{A,P_0}} \leq C^\sharp
$$
for some constants $C_{r_1,r_2}$ and $C^\sharp\geq 1$ independent of
$P_0$. Let
$$
C^\flat= 2n\left({1\over r_1^{2(n+1)}}+1 +C^\sharp\,{1\over
c_A}\,C_{r_1,r_2}\left( {2\,r_2\,C_Y\over
r_2^2-r_1^2}\right)^2\right).
$$
Let $0<r<1$ and let
$$
\hat U_{P_0,r}=U_{P_0}\cap\left\{\left|z^{(P_0)}\right|<{r\over
n\,\sqrt{C^\flat}}\right\}.
$$
Let $N$ be the complex dimension of the subspace of all elements
$$s\in\Gamma\left(Y, L+K_Y+A\right)$$ such that
$$
\int_Y\left|s\right|^2\,e^{-\xi}\,h_A<\infty\; .
$$
Then there exist
$$
s_1,\cdots, s_N\in \Gamma\left(Y, L+K_Y+A\right)
$$
with
$$ \int_Y\left|s_k\right|^2
e^{-\xi}h_A\leq 1
$$
($1\leq k\leq N$) such that, for any $P_0\in Y$ and for any
holomorphic section $s$ of $L+K_Y+A$ over $U_{P_0}$ with
$$
\int_{U_{P_0}}|s|^2 e^{-\xi}h_A=C_s<\infty\; ,
$$
one can find holomorphic functions $b_{P_0,k}$ on $\hat U_{P_0,r}$
such that
$$
s=\sum_{k=1}^Nb_{P_0,k}\,s_k
$$
on $\hat U_{P_0,r}$ and
$$
\sup_{\hat U_{P_0,r}}\sum_{k=1}^N \left|b_{P_0,k}\right|^2 \leq
C^\natural\,C_s\; ,
$$
where $C^\natural={1\over\left(1-r\right)^2}\,C^\flat$.

\medbreak\noindent(Theorem (1.4.4) was proved in [Si02a, p.234, Th.
2.1] for the case where the line bundle $L$ is $mK_Y$. The proof
there works for the case of a general line bundle $L$.  In its
application here and in [Si02a] also, only the case of $L=mK_Y$ is
used.)

\bigbreak\noindent{\sc (1.5) Two-Tower Argument for Invariance of
Plurigenera.}

\bigbreak The detailed proof of the deformational invariance of the
plurigenera for the projective algebraic case was presented in
[Si02a] and [Si03]. The proofs given there are for more general
settings. In order to more transparently present the essence of the
argument, we give here the proof just for the original conjecture
without including any setting which is more general.

\medbreak We now start the proof of the deformational invariance of
the plurigenera for the general case of projective algebraic
manifolds not necessarily of general type.

\medbreak Let $m_0$ be a positive integer. Take
$s^{(m_0)}\in\Gamma\left(X_0,m_0K_{X_0}\right)$. For the proof of
the deformational invariance of the plurigenea we have to extend
$s^{(m_0)}$ to an element of $\Gamma\left(X,m_0K_X\right)$. Let
$A$ be a positive line bundle over $X$ which is sufficiently
positive for the purpose of the global generation of multiplier
ideal sheaves in the sense of (1.4.2) and (1.4.4). Let $h_A$ be a
smooth metric for $A$ with positive curvature form on $X$.  Let
$s_A$ be a global holomorphic section of $A$ over $X$ whose
restriction to $X_0$ is not identically zero.

\medbreak\noindent(1.5.1) {\it Use of Most Singular Metric Still
Giving Finite Norm of Given Initial Section.} Let $\psi={1\over
m_0}\log\left|s^{(m_0)}\right|^2$. Fix arbitrarily
$\ell\in{\mathbb N}$.  Let $m_1=\ell\,m_0$. For $1\leq p<m_1$ let
$\varphi_{p-1}=\left(p-1\right)\psi$ and let $N_p$ be the complex
dimension of the subspace of all elements
$s\in\Gamma\left(X_0,p\,K_{X_0}+A\right)$ such that
$$
\int_{ X_0}\left|s\right|^2\,e^{-\varphi_{p-1}}\,h_A<\infty.
$$
Let $\tilde N=\sup_{1\leq p<m_1}N_p$. Then $\tilde N$ is bounded
independently of $\ell$. This supremum bound being independent of
$\ell$ is a key point in this proof.  (The metric $e^{-\psi}$ is
chosen to be as singular as possible so as to make $N_p$ as small as
possible to guarantee the finiteness of $\tilde N$\,, and yet the
initial section $s^{(m_0)}$ still has finite $L^2$ norm with respect
to the appropriate power of $e^{-\psi}$.)

\medbreak\noindent(1.5.2) {\it Local Trivializations of Line
Bundles.} Let $\{U_\lambda\}_{1\leq\lambda\leq\Lambda}$ be a
finite covering of $X_0$ by Stein open subsets such that

\medbreak\noindent (i) for some $U_\lambda^*$ which contains
$U_\lambda$ as a relatively compact subset, the assumptions and
conclusions of Theorem (1.4.4) are satisfied with each $U_\lambda$
being some $\hat U_{P_0,r}$ from Theorem (1.4.4) when $L$ with
metric $e^{-\xi}$ on $Y$ is replaced by $p\,K_{X_0}$ with metric
$e^{-p\,\psi}$ on $X_0$ for $0\leq p\leq m_1-1$, and

\medbreak\noindent (ii) for $1\leq\lambda\leq\Lambda$ there exists
nowhere zero
$$\xi_\lambda\in\Gamma\left(U_\lambda^{**},-K_{X_0}\right)
$$
for some open subset $U_\lambda^{**}$ of $X_0$ which contains
$U_\lambda^*$ as a relatively compact subset.

\medbreak\noindent(1.5.3) {\it Diagram of ``Two-Tower Argument.''}
First we schematically explain the ``two-tower argument'' and then
we give the details about estimates and convergence.  For the
``two-tower argument'' we start with $\left(s^{(m_0)}\right)^\ell
s_A$ for a large integer $\ell$ at the upper right-hand corner of
the following picture and goes down the tower left of the column
$$
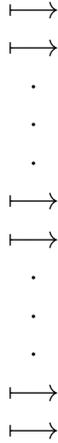
\begin{matrix}\longmapsto\cr\longmapsto\cr\cdot\cr
\cdot\cr\cdot\cr\longmapsto\cr\longmapsto\cr\cdot\cr
\cdot\cr\cdot\cr\longmapsto\cr\longmapsto\cr\end{matrix}
$$
and then goes up the tower right of the column in the following
way.

\bigbreak\noindent
$$
\begin{matrix}\left(s^{(m_0)}\right)^\ell s_A&\quad&\quad&\qquad&&&\widetilde{\sigma_\ell}
\cr\cr \xi_\lambda\left(s^{(m_0)}\right)^\ell
s_A&\quad&s_1^{(m_1-1)},\cdots,s_{N_{m_1-1}}^{(m_1-1)}&\
&\longmapsto&\
&\widetilde{s_1^{(m_1-1)}},\cdots,\widetilde{s_{N_{m_1-1}}^{(m_1-1)}}\cr\cr
\xi_\lambda\,s^{(
m_1-1)}_j&\quad&s_1^{(m_1-2)},\cdots,s_{N_{m_1-2}}^{(m_1-2)}&\
&\longmapsto&\
&\widetilde{s_1^{(m_1-2)}},\cdots,\widetilde{s_{N_{m_1-2}}^{(m_1-2)}}\cr\cr
\cdot&\quad&\cdot&\qquad&&&\cdot\cr
\cdot&\quad&\cdot&\qquad&&&\cdot\cr\cdot&\qquad&\cdot&\qquad&\qquad&&\cdot\cr\cr
\xi_\lambda\,s^{(p+2)}_j&\quad&s_1^{(p+1)},\cdots,s_{N_{p+1}}^{(p+1)}&\
&\longmapsto&\
&\widetilde{s_1^{(p+1)}},\cdots,\widetilde{s_{N_{p+1}}^{(p+1)}}\cr\cr
\xi_\lambda\,s^{(p+1)}_j&\quad&s_1^{(p)},\cdots,s_{N_{p}}^{(p)}&\
&\longmapsto&\
&\widetilde{s_1^{(p)}},\cdots,\widetilde{s_{N_{p}}^{(p)}}\cr\cr
\cdot&\quad&\cdot&\qquad&&&\cdot\cr
\cdot&\quad&\cdot&\qquad&&&\cdot\cr\cdot&\qquad&\cdot&\qquad&&&\cdot\cr\cr
\xi_\lambda\,s^{(3)}_j&\quad&s_1^{(2)},\cdots,s_{N_2}^{(2)}&\
&\longmapsto&\
&\widetilde{s_1^{(2)}},\cdots,\widetilde{s_{N_2}^{(2)}}\cr\cr
\xi_\lambda\,s^{(2)}_j&\quad &s_1^{(1)},\cdots,s_{N_1}^{(1)}&\
&\longmapsto&\
&\widetilde{s_1^{(1)}},\cdots,\widetilde{s_{N_1}^{(1)}}\cr\cr
\end{matrix}
$$

\bigbreak\noindent(1.5.4) {\it Going down the left tower}. We
start out from the top of the left tower with
$$
\left(s^{(m_0)}\right)^\ell \,s_A\in\Gamma\left(X_0,
m_1K_{X_0}+A\right).
$$
We descend one level down the left tower by locally multiplying it
by $\xi_\lambda$ to form $\xi_\lambda\, \left(s^{(m_0)}\right)^\ell
\,s_A$.  By Theorem (1.4.4) (and because each $U_\lambda$ equals
some $\hat U_{P_0,r}$), we can write
$$
\xi_\lambda\, \left(s^{(m_0)}\right)^\ell
\,s_A=\sum_{k=1}^{N_{m_1-1}}
b^{(m_1-1,\lambda)}_k\,s_k^{(m_1-1)}\leqno{(1.5.4.1)}
$$
on $U_\lambda$, where $b^{(m_1-1,\lambda)}_k$ are holomorphic
functions on $U_\lambda$ and where $$
s^{(m_1-1)}_1,\cdots,s^{(m_1-1)}_{N_p}
\in\Gamma\left(X_0,\left(m_1-1\right)K_{X_0}+A\right)
$$
with
$$
\int_{X_0}\left|s^{(m_1-1)}_k\right|^2
e^{-\left(m_1-2\right)\psi}\,h_A\leq 1.
$$

\medbreak We now inductively go down the left tower one level at a
time. For $1\leq p\leq m_1-2$, at the $(p+1)$-st level we have
$$
s^{(p+1)}_j\in\Gamma\left(X_0, (p+1)K_{X_0}+A\right)
$$
and we descend one level down the left tower to the $p$-th level by
locally multiplying it by $\xi_\lambda$ to form $\xi_\lambda\,
s^{(p+1)}_j$.  Again, by Theorem (1.4.4) (and because each
$U_\lambda$ equals some $\hat U_{P_0,r}$), we can write
$$\xi_\lambda\,s^{(p+1)}_j=
\sum_{k=1}^{N_p}b^{(p,\lambda)}_{j,k}s^{(p)}_k\leqno{(1.5.4.2)_p}
$$
on $U_\lambda$, where $ b^{(p,\lambda)}_{j,k}$ are holomorphic
functions on $U_\lambda$ and where $$ s^{(p)}_1,\cdots,s^{(p)}_{N_p}
\in\Gamma\left(X_0,p\,K_{X_0}+A\right)
$$
with
$$
\int_{X_0}\left|s^{(p)}_k\right|^2
e^{-\left(p-1\right)\psi}\,h_A\leq 1.
$$
When we get to the bottom of the left tower, the value of $p$
becomes $1$.

\medbreak Theorem (1.4.4) gives us the estimates:
$$
\displaylines{\sup_{U_\lambda}\sum_{k=1}^{N_{m_1-1}}
\left|b^{(m_1-1,\lambda)}_k\right|^2\leq C^\natural\; ,\cr
\sup_{U_\lambda}\sum_{k=1}^{N_p}\left|b^{(p,\lambda)}_{j,k}\right|^2
\leq C^\natural\cr}
$$
for $1\leq p\leq m_1-2$.  (Here the constant $C^\natural$ is
chosen to have absorbed also the contribution of $\xi$ and the
norm of the initial $\left(s^{(m_0)}\right)^\ell s_A$.) From these
estimates and $(1.5.4.1)$ and $(1.5.4.2)_p$ we obtain the
following estimates:
$$
\int_{X_0}\frac{\left|\left(s^{(m_0)}\right)^\ell
s_A\right|^2}{\max_{1\leq j\leq
N_{m_1-1}}\left|\widetilde{s^{(m_1-1)}_j}\right|^2}\leq C^\natural,
\leqno{(1.5.4.3)}$$
$$\int_{X_0}\frac{\left|s^{(p+1)}_j
\right|^2}{\max_{1\leq j\leq
N_p}\left|\widetilde{s^{(p)}_j}\right|^2}\leq
C^\natural\leqno{(1.5.4.4)_p}
$$
for $1\leq p\leq m_1-2$.

\bigbreak\noindent(1.5.5) {\it Going up the right tower.} Now at the
bottom level we move from the left tower to the right tower and then
move up the right tower one level at a time, by using the extension
theorem of Ohsawa-Takegoshi type (1.4.3).

\medbreak Let $C^\sharp=8\pi\,e\,\sqrt{2+{1\over e}}$.  At bottom
level of $p=1$, because of
$$
\int_{X_0}\left|s^{(1)}_j\right|^2\,h_A\leq 1,
$$
we can extend $$s^{(1)}_j\in\Gamma\left(X_0, K_{X_0}+A\right)$$ to
$$\widetilde{s^{(1)}_j}\in\Gamma\left(X, K_X+A\right)$$ with
$$
\int_X\left|\widetilde{s^{(1)}_j}\right|^2\,h_A\leq C^\sharp.
$$
In the picture of the ``two towers'' the line
$$s_1^{(1)},\cdots,s_{N_1}^{(1)}\ \longmapsto\
\widetilde{s_1^{(1)}},\cdots,\widetilde{s_{N_1}^{(1)}}
$$
signifies the extension of $s^{(1)}_j$ to $\widetilde{s^{(1)}_j}$.

\medbreak Inductively we are going to move up one level at a time
on the right tower.  Suppose we have already moved up to the
$p$-th level for $1\leq p\leq m_1-2$ so that we have the extension
of $$s^{(p)}_j\in\Gamma\left(X_0, pK_{X_0}+A\right)$$ to
$$\widetilde{s^{(p)}_j}\in\Gamma\left(X, pK_X+A\right).$$  We use
$$
\frac{1}{\max_{1\leq j\leq
N_p}\left|\widetilde{s^{(p)}_j}\right|^2}
$$
as the metric for $pK_X+A$ on $X$ and apply the extension theorem of
Ohsawa-Takegoshi type (1.4.3).  From the estimate $(1.5.4.4)_p$ we
can extend
$$s^{(p+1)}_j\in\Gamma\left(X_0, (p+1)K_{X_0}+A\right)$$ to
$$\widetilde{s^{(p+1)}_j}\in\Gamma\left(X, (p+1)K_X+A\right)$$
with $$\int_X{\left|\widetilde{ s^{(p+1)}_j}\right|^2
\over\max_{1\leq k\leq N_p} \left|\widetilde{ s^{(p)}_k}\right|^2}
\leq C^\natural\,C^\sharp.\leqno{(1.5.5.1)_p}$$ In the picture of
the ``two towers'' the line
$$s_1^{(p+1)},\cdots,s_{N_{p+1}}^{(p+1)}\ \longmapsto\
\widetilde{s_1^{(p+1)}},\cdots,\widetilde{s_{N_{p+1}}^{(p+1)}}
$$
signifies the extension of $s^{(p+1)}_j$ to
$\widetilde{s^{(p+1)}_j}$.

\medbreak When we get to the second highest level on the right
tower signified by the line
$$s_1^{(m_1-1)},\cdots,s_{N_{m_1-1}}^{(m_1-1)}\ \longmapsto\
\widetilde{s_1^{(m_1-1)}},\cdots,\widetilde{s_{N_{m_1-1}}^{(m_1-1)}},
$$
we can use
$$
\frac{1}{\max_{1\leq j\leq
N_{m_1-1}}\left|\widetilde{s^{(m_1-1)}_j}\right|^2}
$$
as the metric for $(m_1-1)K_X+A$ on $X$ and apply the extension
theorem of Ohsawa-Takegoshi type (1.4.3).  From the estimate
$(1.5.4.3)$ we can extend
$$\left(s^{(m_0)}\right)^\ell\,s_A\in\Gamma\left(X_0, m_1K_{X_0}+A\right)$$ to
$$\widetilde{\sigma_\ell}\in\Gamma\left(X, m_1K_X+A\right)$$
with
$$\int_X{\left|\widetilde{\sigma_\ell}\right|^2
\over\max_{1\leq k\leq N_{m_1-1}} \left|\widetilde{
s^{(m_1-1)}_k}\right|^2} \leq
C^\natural\,C^\sharp.\leqno{(1.5.5.2)}$$

\medbreak\noindent(1.5.6) {\it Limit Metric and Its Convergence.}
We are going to use
$$
\frac{1}{\limsup_{\ell\to\infty}\left|\widetilde{\sigma_\ell}\right|^{\frac{2}{\ell
m_0}}}
$$
to define a (possibly singular) metric $e^{-\chi}$ for $K_X$ with
$\chi$ plurisubharmonic so that from
$$
\int_{X_0}\left|s^{(m_0)}\right|^2 e^{-(m_0-1)\chi}<\infty,
$$
we get an extension $\widetilde{s^{(m_0)}}\in\Gamma\left(X,
m_0K_X\right)$ of $s^{(m_0)}$.

\bigbreak\noindent{\sc (1.6) Convergence Argument with the Most
Singular Metric.}

\bigbreak We need the convergence of
$$
\limsup_{\ell\to\infty}\left|\widetilde{\sigma_\ell}\right|^{\frac{2}{\ell
m_0}},
$$
which has to come from the $L^2$ estimates $(1.5.5.1)_p$ and
$(1.5.5.2)$. Since the $L^2$ estimate of each factor in a product
(of at least two factors) would not be able to yield any estimate of
the product, we are forced to use the concavity of the logarithmic
function to convert quotients to differences.  For this process it
is essential that the dimension bound $\tilde N$ is independent of
$\ell$. Finally we get pointwise estimates from integral estimates
by using the subharmonicity of the logarithm of the absolute value
of a holomorphic function.  The details of the estimates are as
follows.

\medbreak\noindent(1.6.1) {\it Key Point of Using the Most Singular
Allowable Metric and Uniform Bound of Dimensions of Spaces of
Sections.} For nonnegative valued functions $g_1,\cdots, g_N$, in
general the coefficient $N$ on the right-hand side of the inequality
$$
\int\max\left(g_1,\cdots,g_N\right)\leq N\max\left(\int
g_1,\cdots,\int g_N\right)
$$
cannot be lowered, because the supports of $g_1,\cdots,g_N$ may all
be disjoint or close to being disjoint.  So from $(1.5.5.1)_p$ we
get
$$\int_X{\max_{1\leq j\leq N_{p+1}}\left|\widetilde{ s^{(p+1)}_j}\right|^2
\over\max_{1\leq k\leq N_p} \left|\widetilde{ s^{(p)}_k}\right|^2}
\leq N_{p+1}\, C^\natural\,C^\sharp,
$$
where the factors $N_{p+1}$ on the right-hand side cannot be lowered
in general.  Since $\tilde N=\sup_{1\leq p<m_1}N_p$, it follows that
$$
\int_X{\max_{1\leq j\leq N_{p+1}}\left|\widetilde{
s^{(p+1)}_j}\right|^2 \over\max_{1\leq k\leq N_p} \left|\widetilde{
s^{(p)}_k}\right|^2} \leq\tilde N\,
C^\natural\,C^\sharp.\leqno{(1.6.1.1)_p}
$$
The main point is that the right-hand side $\tilde N\,
C^\natural\,C^\sharp$ is now independent of $p$, which is essential
for the estimate in passing to the final limit metric.

\medbreak\noindent(1.6.2) {\it Preliminary Notations for Local
Trivializations of Line Bundles.} Let $0<r_0<r_1<r_2<1$. Choose a
finite number of coordinate charts $\tilde W_\lambda$ in $X$ with
coordinates
$$
\left(z^{(\lambda)},t\right)=
\left(z^{(\lambda)}_1,\cdots,z^{(\lambda)}_n,t\right)
$$
for $1\leq\lambda\leq\hat\Lambda$ such that

\medbreak\noindent (ii) each $$W_\lambda^{\prime\prime}:=
\left\{\left|z^{(\lambda)}_1\right|<r_2,
\cdots,\left|z^{(\lambda)}_n\right|<r_2,|t|<r_2 \right\}
$$
is relatively compact in $\tilde W_\lambda$ for
$1\leq\lambda\leq\hat\Lambda$,

\medbreak\noindent (ii)
$X\cap\left\{\left|t\right|<r_2\right\}=\bigcup_{\lambda=1}^{\hat\Lambda}
W_\lambda^{\prime\prime}$,

\medbreak\noindent (iii)
$X\cap\left\{\left|t\right|<r_0\right\}=\bigcup_{\lambda=1}^{\hat\Lambda}
W_\lambda$, where
$$W_\lambda=
\left\{\left|z^{(\lambda)}_1\right|<r_0,
\cdots,\left|z^{(\lambda)}_n\right|<r_0,|t|<r_0\right\}
$$
for $1\leq\lambda\leq\hat\Lambda$, and

\medbreak\noindent (iv) there exist nowhere zero
$$
\displaylines{
$$
\hat\tau_{\lambda,A}\in\Gamma\left(W_\lambda^{\prime\prime},A\right),\cr
\hat\xi_\lambda\in\Gamma\left(W_\lambda^{\prime\prime},-K_X\right)\cr
}
$$
for $1\leq\lambda\leq\hat\Lambda$.

\medbreak\noindent Let $dV_{z^{(\lambda)},t}$ be the Euclidean
volume form in the coordinates system
$\left(z^{(\lambda)},t\right)$.  Let
$$W_\lambda^\prime=
\left\{\left|z^{(\lambda)}_1\right|<r_1,
\cdots,\left|z^{(\lambda)}_n\right|<r_1,|t|<r_1 \right\}.
$$

\medbreak\noindent(1.6.3) {\it Concavity of Logarithm for
Conversion from Integrals of Quotients to Differences of
Integrals.} From the concavity of the logarithmic function we
conclude that
$$
\displaylines{{1\over\left(\pi r_1^2 \right)^{n+1}}\,\int_{
W_\lambda^\prime}\left(\log\max_{1\leq j\leq N_{p+1}}\left|
\hat\tau_{\lambda,A}\,\hat\xi_\lambda^{p+1}\,\tilde
s^{(p+1)}_j\right|^2\right)dV_{z^{(\lambda)},t} \cr -{1\over
\left(\pi r_1^2\right)^{n+1}}\,\int_{W_\lambda^\prime}\left(\log
\max_{1\leq k \leq N_p}\left|
\hat\tau_{\lambda,A}\,\hat\xi_\lambda^p\,\tilde
s^{(p)}_k\right|^2\right)dV_{z^{(\lambda)},t}\cr = {1\over \left(\pi
r_1^2\right)^{n+1}}\,\int_{W_\lambda^\prime}\left(\log{\max_{1\leq
j\leq N_{p+1}}\left|
\hat\tau_{\lambda,A}\,\hat\xi_\lambda^{p+1}\,\tilde
s^{(p+1)}_j\right|^2\over \max_{1\leq k \leq N_p}\left|
\hat\tau_{\lambda,A}\,\hat\xi_\lambda^p\,\tilde
s^{(p)}_k\right|^2}\right)dV_{z^{(\lambda)},t}\cr \leq\log
\left({1\over\left(\pi r_1^2
\right)^{n+1}}\,\int_{W_\lambda^\prime}{\max_{1\leq j\leq
N_{p+1}}\left| \hat\tau_{\lambda,A}\,\hat\xi_\lambda^{p+1}\,\tilde
s^{(p+1)}_j\right|^2\over \max_{1\leq k \leq N_p}\left|
\hat\tau_{\lambda,A}\,\hat\xi_\lambda^p\,\tilde
s^{(p)}_k\right|^2}\,dV_{z^{(\lambda)},t}\right)\cr
\leq\log\left(\sup_{W_\lambda^\prime} \left({1\over \left(\pi
r_1^2\right)^{n+1}}\left|
\,\hat\xi_\lambda\right|^2\,dV_{z^{(\lambda)},t}\right)
\int_{W_\lambda^\prime}{\max_{1\leq j\leq N_{p+1}}\left|\, \tilde
s^{(p+1)}_j\right|^2\over \max_{1\leq k \leq N_p}\left| \tilde
s^{(p)}_k\right|^2}\right)\cr \leq \log\left(\tilde
N\,C^\natural\,C^\sharp\,\sup_{W_\lambda^\prime}\left({1\over
\left(\pi r_1^2\right)^{n+1}}\left|
\,\hat\xi_\lambda\right|^2\,dV_{z^{(\lambda)},t} \right)\right)\cr }
$$
for $1\leq p\leq m_1-2$, where $(1.6.1.1)_p$ is used for the last
inequality. Likewise
$$
\displaylines{ {1\over\left(\pi r_1^2
\right)^{n+1}}\,\int_{W_\lambda^\prime}\left(\log\left|
\hat\tau_{\lambda,A}\,\hat\xi_\lambda^{m_1}\,\tilde\sigma_\ell\right|^2\right)dV_{z^{(\lambda)},t}
\cr -{1\over\left(\pi r_1^2
\right)^{n+1}}\,\int_{W_\lambda^\prime}\left( \log\max_{1\leq k \leq
N_{m_1-1}}\left|
\hat\tau_{\lambda,A}\,\hat\xi_\lambda^{m_1-1}\,\tilde
s^{(m_1-1)}_k\right|^2\right) dV_{z^{(\lambda)},t}\cr \leq
\log\left(C^\natural\,C^\sharp\,\sup_{W_\lambda^\prime}\left({1\over\left(\pi
r_1^2 \right)^{n+1}}\left|
\,\hat\xi_\lambda\right|^2\,dV_{z^{(\lambda)},t} \right)\right)\cr}
$$
Adding up, we get
$$
\displaylines{ {1\over\left(\pi
r_1^2\right)^{n+1}}\,\int_{W_\lambda^\prime}\left(\log\left|
\hat\tau_{\lambda,A}\,\hat\xi_\lambda^{m_1}\,\tilde\sigma_\ell\right|^2\right)dV_{z^{(\lambda)},t}
\cr \leq {1\over\left(\pi
r_1^2\right)^{n+1}}\,\int_{W_\lambda^\prime}\left( \log\max_{1\leq k
\leq N_1}\left| \hat\tau_{\lambda,A}\,\hat\xi_\lambda\,\tilde
s^{(1)}_k\right|^2\right) dV_{z^{(\lambda)},t}\cr +
\left(m_1-1\right)\log\left(\tilde
N\,C^\natural\,C^\sharp\,\sup_{W_\lambda^\prime}\left({1\over\left(\pi
r_1^2\right)^{n+1}}\left|
\,\hat\xi_\lambda\right|^2\,dV_{z^{(\lambda)},t} \right)\right).\cr}
$$

\medbreak\noindent(1.6.4) {\it Use of Sub-Mean-Value Property of
Logarithm of Absolute-Value of Holomorphic Functions.} By the
sub-mean-value property of plurisubharmonic functions
$$
\displaylines{
\sup_{W_\lambda}\log\left|\hat\tau_{\lambda,A}\,\hat\xi_\lambda^{m_1}\,\tilde\sigma_\ell\right|^2\cr
\leq{1\over\left(\pi\left(r_1-r_0\right)^2\right)^{n+1}}\,
\int_{W_\lambda^\prime}\left(\log\left|\hat\tau_{\lambda,A}\,\hat\xi_\lambda^{m_1}\,\tilde\sigma_\ell\right|^2\right)\,dV_{z^{(\lambda)},t}.\cr
}
$$
Choose a positive number $C^\spadesuit$ such that ${1\over m_0}\,
\log\,C^\spadesuit$ is no less than
$$
\left({r_1\over r_1-r_0}\right)^{2(n+1)}\,\log\left(\tilde
N\,C^\natural\,C^\sharp\, \sup_{
W_\lambda^\prime}\left({1\over\left(\pi r_1^2\right)^{n+1}}\left|
\,\hat\xi_\lambda\right|^2 \,dV_{z^{(\lambda)},t} \right)\right)
$$
for every $1\leq\lambda\leq\Lambda$ and $1\leq p\leq m_1-1$. Since
$\tilde N$ is bounded independently of $\ell$, we can choose
$C^\spadesuit$ to be independent of $\ell$. Let $\hat C$ be
defined by
$$
\log\,\hat C = {1\over\left(\pi r_1\right)
^{n+1}}\,\int_{W_\lambda^\prime}\left( \log\max_{1\leq k \leq
N_1}\left| \hat\tau_{\lambda,A}\,\hat\xi_\lambda\,\tilde
s^{(1)}_k\right|^2\right) dV_{z^{(\lambda)},t}\,.
$$
Then
$$
\sup_{1\leq\lambda\leq\hat\Lambda}
\sup_{W_\lambda}\left|\hat\tau_{\lambda,A}
\,\hat\xi_\lambda^{\ell\,m_0}\,\tilde\sigma_\ell\right|^2 \leq\hat
C\,\left(C^\spadesuit\right)^{\ell\,m_0-1}.
$$

\medbreak\noindent(1.6.5) {\it Final Step of Construction of Limit
Metric.} For $1\leq\lambda\leq\hat\Lambda$ let $\chi_\lambda$ be
the function on $W_\lambda$ which is the upper semi-continuous
envelope of
$$
\limsup_{\ell\rightarrow\infty}\log\left|
\,\hat\xi_\lambda^{\ell\,m_0}\,\tilde\sigma_\ell\right|^{2\over\ell}
$$
which is $\leq m_0 \log\,C^\spadesuit$. From the definition of
$\chi_\lambda$ and the fact that $\tilde\sigma_\ell$ is the
extension of $\left(s^{(m_0)}\right)^\ell s_A$, we have
$$
\sup_{X_0\cap W_\lambda}\left(\left|
\,\hat\xi_\lambda^{m_0}\,s^{(m_0)}\right|^2\,e^{-\chi_\lambda}\right)
\leq 1
$$
for $1\leq\lambda\leq\hat\Lambda$. Let
$e^{-\chi}=\left|\hat\xi_\lambda\right|^2e^{-\frac{\chi_\lambda}{m_0}}$
be the metric of $K_X$ on $X\cap\left\{\left|t\right|<r\right\}$ so
that the square of the pointwise norm of a local section $\sigma$ of
$K_X$ on $W_\lambda$ is
$\left|\sigma\hat\xi_\lambda\right|^2e^{-\chi_\lambda}$.

\medbreak Let
$\left\{\rho_\lambda\right\}_{1\leq\lambda\leq\tilde\Lambda}$ be a
partition of unity subordinate to the open cover $\left\{X_0\cap
W_\lambda\right\}_{1\leq\lambda\leq\tilde\Lambda}$ of $X_0$. Since
$$
\displaylines{\left|s^{(m_0)}\right|^2\,e^{-\left(m_0-1\right)\chi}
=\left|s^{(m_0)}\right|^2\,\left|\tilde\xi_\lambda\right|^{2\left(m_0-1\right)}\,e^{-\left(\frac{m_0-1}{m_0}\right)\chi_\lambda}\cr
=\left|s^{(m_0)}\,
\tilde\xi_\lambda^{m_0}\right|^2\left(\frac{e^{-\left(\frac{m_0-1}{m_0}\right)\chi_\lambda}}{\left|\tilde\xi_\lambda\right|^2}\right)
\leq\frac{e^{-\left(\frac{m_0-1}{m_0}\right)\chi_\lambda}}{\left|\tilde\xi_\lambda\right|^2}\cr}
$$
on $X_0\cap W_\lambda$ for $1\leq\lambda\leq\tilde\Lambda$, it
follows that
$$
\displaylines{\int_{X_0}\left|s^{(m_0)}\right|^2\,e^{-\left(m_0-1\right)\chi}=\sum_{\lambda=1}^{\tilde\Lambda}\int_{X_0\cap
W_\lambda}\rho_\lambda\left|s^{(m_0)}\right|^2\,e^{-\left(m_0-1\right)\chi}\cr\leq\sum_{\lambda=1}^{\tilde\Lambda}\int_{X_0\cap
W_\lambda}\rho_\lambda\left(\frac{e^{-\left(\frac{m_0-1}{m_0}\right)\chi_\lambda}}{\left|\tilde\xi_\lambda\right|^2}\right)<\infty.\cr}
$$
Since the curvature current of the metric $e^{-\chi}$ of $K_X$ on
$X\cap\left\{\left|t\right|<r\right\}$ is nonnegative, by the
extension theorem of Ohsawa-Takegoshi type (1.4.3), we can extend
$s^{(m_0)}$ to an element of
$\Gamma\left(X\cap\left\{\left|t\right|<r\right\}, m_0K_X\right)$.
Finally, by the coherence of the zeroth direct image of ${\mathcal
O}_X\left(m_0K_X\right)$ under the projection $\pi:X\to\Delta$ and
by the Steinness of $\Delta$, we can extend $s^{(m_0)}$ to an
element of $\Gamma\left(X\cap\left\{\left|t\right|<r\right\},
m_0K_X\right)$.  This concludes the proof of the deformational
invariance of the plurigenera for the general case of projective
algebraic manifolds not necessarily of general type.

\bigbreak\noindent{\sc (1.7) Maximally Regular Metric for the
Canonical Line Bundle.}

\bigbreak As we have seen in the above proof, the main idea of the
proof of the deformational invariant of the plurigenera is to
produce a metric $e^{-\chi}$ for $K_X$ with $\chi$ plurisubharmonic
so that a given element $s^{(m_0)}\in\Gamma\left(X_0, m_0
K_X\right)$ satisfies
$$
\int_{X_0}\left|s^{(m_0)}\right|^2\,e^{-\left(m_0-1\right)\chi}<\infty.\leqno{(1.7.0.1)}
$$
The metric $e^{-\chi}$ is constructed as
$$
\frac{1}{\sum_j\left|\tau_j\right|^2},
$$
where each $\tau_j$ is a multi-valued holomorphic section of $K_X$
(in the sense that $\left(\tau_j\right)^q$ for some appropriate
positive integer $q$ is a global holomorphic section of $qK_X$ over
$X$). The metric $e^{-\chi}$ is singular at points of common zeroes
of $\tau_j$.  In order for $(1.7.0.1)$ to be satisfied, it is
natural to choose $\tau_j$ so that $\sum_j\left|\tau_j\right|^2$ is
as large as possible, or equivalently, the metric $e^{-\chi}$ is as
regular as possible.

\medbreak A maximally regular metric for the canonical line bundle
$K_{X_0}$ of $X_0$ is given by $\Phi^{-1}$ with
$$
\Phi=\sum_{m=1}^\infty\varepsilon_m\left(\sum_{j=1}^{q_m}\left|\sigma^{(m)}_j\right|^2\right)^{\frac{1}{m}},
$$
where
$$
\sigma^{(m)}_1,\cdots,\sigma^{(m)}_{q_m}\in\Gamma\left(X_0,
mK_{X_0}\right)
$$
form a basis over ${\mathbb C}$ and $\varepsilon_m$ is a sequence
of positive numbers which decrease to $0$ fast enough to make the
above defining series for $\Phi$ converge.  With the maximally
regular metric $\Phi^{-1}$ for the canonical line bundle $K_{X_0}$
of $X_0$ the condition $(1.7.0.1)$ is automatically satisfied. The
metric $\Phi^{-1}$ is only for $K_{X_0}$ and not for $K_X$. With
the introduction of a holomorphic line bundle $A$ on $X$ which is
sufficiently positive for the global generation of multiplier
ideal sheaves (1.4.2), the ``two-tower argument'' above could give
us
$$
\hat s^{(m)}_1,\cdots,\hat s^{(m)}_{\hat q_m}\in\Gamma\left(X,
mK_X+A\right)
$$
whose restrictions to $X_0$ would give a basis of $\Gamma\left(X_0,
mK_{X_0}+A\right)$ over ${\mathbb C}$.  We could form the metric
$\left(\hat\Phi_m\right)^{-1}$ for $mK+A$ on $X$ with
$$\hat\Phi_m=\sum_{j=1}^{\hat
q_m}\left|\hat s^{(m)}_j\right|^2.
$$
The condition
$$
\int_{X_0}\left|s^{(m_0)}s_A\right|^2\,\left(\hat\Phi_{m_0-1}\right)^{-1}h_A<\infty
$$
is clearly satisfied.

\medbreak A natural way to get rid of the undesirable summand $A$ is
to use the limit of
$$
\left(\hat\Phi_{\left(m_0-1\right)\ell}\right)^{-\,\frac{1}{\ell}}
$$
as a metric of $\left(m_0-1\right)K_X$ when some sequence $\ell$ of
positive integers goes to infinity.  The main obstacle is the
convergence of the limit. The extension of orthonormal basis from
the initial fibers to the whole family in general would not preserve
the orthonormality property.  The use of sub-mean-value property of
the logarithm of the absolute value of holomorphic functions to get
uniform bounds would involve shrinking the domain, which might
eventually disappear completely in an infinite inductive process.
The known techniques of estimates could not yield the needed
convergence.

\medbreak A special method was introduced in [Siu01] to solve this
problem in the special case when manifolds $X_t$ in the family
$\pi:X\to\Delta$ are of general type. For the special case of
general type we can write $aK_X=A+D$ for some positive integer $a$
and some effective divisor $D$ of $X$ not containing any fiber
$X_t$ (after replacing $\Delta$ by a smaller disk if necessary).
Let $s_D$ be the canonical section of the divisor $D$.  We can use
$$
\left(\hat
\Phi_{m_0\ell}\left|s_D\right|^2\right)^{-\,\frac{m_0-1}{m_0\ell+a}}
$$
as the metric for $\left(m_0-1\right)K_X$ for $\ell$ sufficiently
large without passing to the limit with $\ell\to\infty$.   When
$\ell$ is sufficiently large, the condition $(1.7.0.1)$ follows from
using H\"older's inequality to separate the factor
$\left(\left|s_D\right|^2\right)^{-\,\frac{m_0-1}{m_0\ell+a}}$ into
an integral of a power of it.  The use of this particular step
circumvents the obstacle of proving the convergence of the limit.

\medbreak Finally a solution of the convergence problem came in
[Si02a] when it was realized that in the construction of the metric
for $\left(m_0-1\right)K_X$ the priority should be given to the
convergence condition instead of to the finiteness condition
$(1.7.0.1)$.  Instead of using the most natural maximally regular
metric for the canonical line bundle, in [Si02a] the most singular
metric for the canonical line bundle is chosen which still fulfils
the condition $(1.7.0.1)$.

\medbreak Though surprisingly it turns out that the most singular
metric should be used for the deformational invariance of the
plurigenera, the maximally regular metric of the canonical line
bundle will play an important r\^ole in the problem of the finite
generation of the canonical ring which we will discuss below.

\bigbreak\noindent{\sc (1.8) Finite Generation of Canonical Ring
and Skoda's Estimates for Ideal Generation.}

\bigbreak We now turn to the problem of the finite generation of the
canonical ring.  A very powerful tool for this problem is the
following estimate of Skoda for ideal generation [Sk72].

\medbreak\noindent(1.8.1) {\it Theorem (Skoda)}. Let $\Omega$ be a
pseudoconvex domain spread over ${\mathbb C}^n$ and $\psi$ be a
plurisubharmonic function on $\Omega$.  Let $g_1,\ldots,g_p$ be
holomorphic functions on $\Omega$.  Let $\alpha>1$ and
$q=\inf(n,p-1)$.  Then for every holomorphic function $f$ on
$\Omega$ such that
$$
\int_\Omega|f|^2|g|^{-2\alpha q-2}e^{-\psi}d\lambda<\infty,
$$
there exist holomorphic functions $h_1,\ldots,h_p$ on $\Omega$
such that
$$
f=\sum_{j=1}^pg_jh_j
$$
and
$$
\int_\Omega|h|^2|g|^{-2\alpha q}e^{-\psi}d\lambda
\leq\frac{\alpha}{{\alpha-1}}\int_\Omega|f|^2|g|^{-2\alpha
q-2}e^{-\psi}d\lambda,
$$
where
$$
|g|=\left(\sum_{j=1}^p|g_j|^2\right)^{\frac{1}{2}},\qquad
|h|=\left(\sum_{j=1}^p|h_j|^2\right)^{\frac{1}{2}},
$$
and $d\lambda$ is the Euclidean volume element of ${\mathbb C}^n$.

\bigbreak\noindent(1.8.2) {\it Skoda's Estimate Applied to Sections
of Line Bundles.}  By using the procedure which reduces the
situation of sections of line bundles over compact algebraic
manifolds to functions on Stein domains over ${\mathbb C}^n$ with
$L^2$ estimates, we can translate Skoda's estimate into the
following algebraic geometric formulation.

\medbreak\noindent(1.8.3) {\it Theorem.} Let $X$ be a compact
complex algebraic manifold of dimension $n$, $L$ a holomorphic
line bundle over $X$, and $E$ a holomorphic line bundle on $X$
with metric $e^{-\varphi }$ ($\varphi$ possibly asssuming
$-\infty$ values) such that $\psi$ is plurisubharmonic. Let $k\geq
1$ be an integer, $G_{1},\ldots ,G_{p}\in \Gamma (X,L)$, and
$\left\vert G\right\vert ^{2}=\underset{j=1}{\overset{p}{\sum
}}\left\vert G_{j}\right\vert ^{2}$. Let ${\mathcal I}={\mathcal
I}_{(n+k+1)\log\left|G\right|^2+\varphi}$ and ${\mathcal
J}={\mathcal I}_{(n+k)\log\left|G\right|^2+\varphi }$ Then
$$\Gamma \left(X,{\mathcal I}\otimes\left((n+k+1)L+E+K_X\right)\right)
=\sum_{j=1}^pG_{j}\Gamma \left(X,{\mathcal J}\otimes\left(
(n+k)L+E+K_X\right)\right).$$

\medskip
\noindent {\it Proof}. Clearly the right-hand side is contained in
the left-hand side of the equality in the conclusion.  To prove the
opposite direction, we take $F\in\Gamma (X,\mathcal{I}\otimes
(n+k+1)L+E+K_{X})$.

\medbreak Let $S$ be a non identically zero meromorphic section of
$E$ on $X$.  We take a branched cover map $\pi:X\rightarrow{\bf
P}_n$. Let $Z_0$ be a hypersurface in ${\bf P}_n$ which contains the
infinity hyperplane of ${\bf P}_n$ and the branching locus of $\pi$
in ${\bf P}_n$ such that $Z:=\pi^{-1}(Z_0)$ contains the divisor of
$G_1$ and the zero-set of $S$ and the pole-set of $S$. Let
$\Omega=X-Z$.

\medbreak Let $$g_j=\frac{G_j}{G_1}\qquad(1\leq j\leq p)$$ and
define $f$ by
$$
\frac{F}{G_1^{n+k+1}S}=fdz_1\wedge\cdots\wedge dz_n,
$$
where $z_1,\cdots,z_n$ are the affine coordinates of ${\mathbb
C}^n$.  Use $\alpha=\frac{n+k}{n}$.  Let $\psi=\varphi-\log|S|^2$.
It follows from $F\in{\cal I}_{(n+k+1)\log|G|^2+\varphi}$ that
$$
\int_X\frac{|F|^2}{|G|^{2(n+k+1)}}e^{-\varphi}<\infty,
$$
which implies that
$$
\int_\Omega\frac{|f|^2}{|g|^{2(n+k+1)}}e^{-\psi}=
\int_\Omega\frac{\left|\frac{F}{G_1^{n+k+1}S}\right|^2}{\left|\frac{G}{
G_1}\right|^{2(n+k+1)}}e^{-\psi}=
\int_\Omega\frac{|F|^2}{|G|^{2(n+k+1)}}e^{-\varphi}<\infty,
$$
where $\left|g\right|^2=\sum_{j=1}^p\left|g_j\right|^2$. By Theorem
(1.8.1) with $q=n$ (which we assume by adding some
$$F_{p+1}\equiv\cdots\equiv F_{n+1}\equiv 0$$ if $p<n+1$) so that
$$2\alpha q+2=2\cdot\frac{n+k}{n}\cdot n+2=2(n+k+1).$$  Thus there
exist holomorphic functions $h_1,\cdots,h_p$ on $\Omega$ such that
$f=\sum_{j=1}^p g_jh_j$ and
$$
\sum_{j=1}^p\int_\Omega\frac{|h_j|^2}{|g|^{2(n+k)}}e^{-\psi}<\infty.
$$
Define
$$
H_j=G_1^{n+k}h_j S dz_1\wedge\cdots\wedge dz_n.
$$
Then
$$
\int_\Omega\frac{|H_j|}{|G|^{2(n+k)}}e^{-\varphi}=
\int_\Omega\frac{|h_j|}{|g|^{2(n+k)}}e^{-\psi}<\infty
$$
so that $H_j$ can be extended to an element of $\Gamma
\left(X,{\mathcal J}\otimes\left( (n+k)L+E+K_X\right)\right)$.
Q.E.D.

\bigbreak As an illustration of how Skoda's estimate could be used
to yield readily results on finite generation, we give the following
statement in the case of globally free line bundles.

\medbreak\noindent(1.8.4) {\it Theorem (Generation for the Ring of
Sections of Multiples of Free Bundle).} Let $F$ be a holomorphic
line bundle over a compact projective algebraic manifold $X$ of
complex dimension $n$. Let $a>1$ and $b\geq 0$ be integers such that
$aF$ and $bF-K_X$ are globally free over $X$. Then the ring
$\oplus_{m=0}^\infty\Gamma(X,mF)$ is generated by
$\oplus_{m=0}^{(n+2)a+b-1}\Gamma(X,mF)$.

\medskip
\noindent {\it Proof.} For $0\leq\ell<a$ let
$E_\ell=(b+\ell)F-K_X$ and $L=aF$.  Let $G_1,\cdots,G_p$ be a
basis of $\Gamma(X,L)=\Gamma(X,aF)$.  Let $H_1,\cdots,H_q$ be a
basis of $\Gamma(X,bF-K_X)$.  We give $E_\ell$ the metric
$$
\frac{1}{(\sum_{j=1}^p|G_j|^{\frac{2\ell}{
a}})(\sum_{j=1}^q|H_j|^2)}.
$$
Since both ${\mathcal I}$ and ${\mathcal J}$ (from Theorem (1.8.4)
when $E$ is set to be $E_\ell$) are unit ideal sheaves due to the
global freeness of $aF$ and $bF-K_X$, it follows from Theorem
(1.8.4) that
$$
\Gamma(X,(n+k+1)L+E_\ell+K_X) =\sum_{j=1}^p
G_j\Gamma(X,(n+k)L+E+K_X)
$$
for $k\geq 1$ and $0\leq\ell<a$, which means that
$$
\Gamma\left(X,((n+k+1)a+\ell+b)F\right) =\sum_{j=1}^p
G_j\Gamma\left(X,((n+k)a+\ell+b)F\right)
$$
for $k\geq 1$ and $0\leq\ell<a$.  Thus
$\oplus_{m=0}^{(n+2)a+b-1}\Gamma(X,mF)$ generates the ring
$\oplus_{m=0}^\infty\Gamma(X,mF)$.  Q.E.D.

\bigbreak\noindent{\sc (1.9) Stable Vanishing Order.}  The use of
Skoda's esimates in the simple case of finite generation for the
ring of sections of multiples of free line bundles suggests that the
finite generation of the canonical ring for the case of finite type
depends on properties of the stable vanishing order which is defined
as follows.

\medbreak Let
$$
s^{(m)}_1, \cdots, s^{(m)}_{q_m}\in\Gamma\left(X,mK_X\right)
$$
be a basis over ${\mathbb C}$.  Let $$\Phi =\sum_{m=1}^\infty
\varepsilon _{m}\left( \sum\limits_{j=1}^{q_{m}}\left\vert
s_{j}^{(m)}\right\vert ^{2}\right) ^{\frac{1}{m}},$$ where
$\varepsilon_m$ is a sequence of positive numbers decreasing fast
enough to guarantee convergence of the series in the definition of
$\Phi$.

\medbreak We would like to remark that $\frac{1}{\Phi }$ is the
metric for $K_{X}$ which was introduced earlier in our discussion
of the deformational invariance of plurigenera.

\bigbreak One main difficulty of applying Skoda's estimates to
prove the conjecture on the finite generation of the canonical
ring is the common zero-set of $\Phi$ and the vanishing orders of
$\Phi$ at its points.  A more precise definition of the vanishing
order of $\Phi$ is the Lelong number of the closed positive
$(1,1)$-current
$$
\frac{\sqrt{-1}}{2\pi}\,\partial\bar\partial\log\Phi
$$
(see [Si74] for the definition and properties of Lelong numbers). We
call the vanishing orders of $\Phi$ the {\it stable vanishing
orders.} When the canonical ring is finitely generated, it is clear
that all vanishing orders of $\Phi$ are rational.  As an
intermediate step of the proof of the finite generation of the
canonical ring, one has the following conjecture on the rationality
of the stable vanishing orders.

\bigbreak\noindent(1.9.1) {\it Conjecture on Rationality of Stable
Vanishing Orders.} If $X$ is of general type, then all the
vanishing orders of $\Phi$ are rational.

\bigbreak\noindent(1.9.2) {\it Approach to Rationality by Degenerate
Complex Monge-Amp\`ere Equation.}  One approach to the conjecture on
the rationality of stable vanishing orders is by degenerate complex
Monge-Amp\`ere equations and indicial equations.  Let $\Omega$ be
the open subset of $X$ consisting of all points where $\Phi$ is
positive. The degenerate complex Monge-Amp\`ere equation is
$$
(\sqrt{-1}\,\partial\overline{\partial }\log \ \omega )^{n}=\omega
\text{ \ on }\Omega\leqno{(1.9.2.1)}$$ with the boundary condition
that the quotient
$$\frac{\left(\sqrt{-1}\partial\bar\partial\Phi\right)^n}{\omega}$$
is bounded from above and below by positive constants on
$\Omega\cap U$ for some neighborhood $U$ of $\partial\Omega$,
where the unknown $\omega$ is a positive smooth $(n,n)$-form on
$\Omega$ so that $\sqrt{-1}\,\partial\bar\partial\log\omega$ is
strictly positive on $\Omega$.

\medbreak The motivation of considering the degenerate complex
Monge-Amp\'ere equation comes from the following two considerations.

\medbreak\noindent(i) According to an observation of Demailly,
another way of getting a metric of $K_X$ comparable to $\Phi^{-1}$
is consider the maximum $\varphi$ among all metrics $e^{-\varphi}$
of $K_X$ with $\varphi$ locally plurisubharmonic so that
$e^{-\varphi}h$ has infimum $1$ on $X$ for some smooth metric $h$
of $K_X$.  The reason is that one can use $L^2$ estimates of
$\bar\partial$ to get global holomorphic sections of a line bundle
$L$ with semipositive curvature current after twisting $L$ by a
sufficiently positive line bundle $A$ independent of $L$.
Conversely, global sections can be used to define a metric with
semipositive curvature current.

\medbreak\noindent(i) By the work of Bedford and Taylor on the
Dirichlet problem for complex Monge-Amp\`ere equations [BT76], a
solution of the complex Monge-Amp\`ere equation can be obtained by
the Perron method of maximization subject to certain normalization.

\medbreak We would like to remark that our degenerate complex
Monge-Amp\`ere equation (1.9.2.1) is very different from the kind of
degenerate complex Monge-Amp\`ere equations considered by Yau in the
last part of his paper [Ya78].

\bigbreak\noindent(1.9.3) {\it Indicial Equations of Regular
Singular Ordinary Differential Equations.}

\bigbreak In order to use the degenerate complex Monge-Amp\`ere
equation to conclude the rationality of the stable vanishing orders,
one compares the situation to the method of using undetermined
coefficients to get vanishing orders of solutions of regular
singular ordinary differential equation.

\medbreak In the case of the regular singular ordinary differential
equation
$$x^{2}y^{\prime \prime }+xa(x)y^{\prime }+b(x)y=0\leqno{(1.9.3.1)}$$
of a single real dependent variable $y$ and a single real
independent variable $x$, for a solution of the form
$$y=x^r\sum_{j=0}^\infty c_jx^j$$
with undetermined coefficients $c_j$ to satisfy (1.9.3.1) the
exponent $r$, known as the {\it index}, has to satisfy the {\it
indicial equation} $$r(r-1)+ra(0)+b(0)=0.\leqno{(1.9.3.2)}$$ The
indicial equation is {\it quadratic} because in general there are
{\it two} solutions of the second-order ordinary equation.

\bigbreak\noindent(1.9.4) {\it Analog of Indicial Equation for
Degenerate Complex Monge-Amp\`ere Equation.}

\bigbreak For our {\it degemerate} complex Monge-Amp\`ere equation
(1.9.2.1), the analog of the indicial equation (1.9.3.2) is a system
of equations whose unknowns are the stable vanishing orders.

\bigbreak Because the solution of the degenerate complex
Monge-Amp\`ere equation is expected to be unique, instead of the
quadratic indical equation (1.9.3.2) from the ordinary differential
equation, the system of equations for the stable vanishing orders is
expected to be linear with rational coefficients, independent enough
to give the rationality of the stable vanishing orders.
\bigbreak\noindent

\bigbreak\noindent(1.9.5) {\it Finite Generation of Canonical Ring
with Assumption of Rationality of Stable Vanishing Orders.}

\bigbreak Suppose we are able to show that all stable vanishing
orders are rational.  How far away are we from concluding the finite
generation of the canonical ring?  To answer this question, we
introduce first the notion of Lelong sets.

\medbreak For a positive number $c$, the Lelong set
$E_c\left(\Phi\right)$ of $\Phi$ (or more precisely,
$\frac{\sqrt{-1}}{2\pi}\,\partial\bar\partial\log\Phi$) means the
set of all points $x$ of $X$ such that the Lelong number of the
closed positive $(1,1)$-current $
\frac{\sqrt{-1}}{2\pi}\,\partial\bar\partial\log\Phi$ is $\geq c$ at
the point $x$.  Note that $E_c(\Phi)$ is always a complex-analytic
subvariety [Si74].

\medbreak An irreducible Lelong set $E$ of $\Phi$ means a branch of
$E_c(\Phi)$ for some $c>0$.  The generic Lelong number on a Lelong
set $E$ is the Lelong number at a generic point of $E$ (which is
independent of the choice of the generic point of $E$).  The
following statement holds.

\bigbreak\noindent(1.9.6) {\it Statement.} If the number of all
irreducible Lelong sets of $\Phi$ is finite and if their generic
Lelong numbers are all rational, then the canonical ring is finitely
generated.

\bigbreak So after the verification of the rationality of all
stable vanishing orders, what remains is the problem of handling
an infinite number of irreducible Lelong sets.  For that one uses
the fact that an additional copy of the canonical line bundle is
added to play the r\^ole of the volume form when one considers
vanishing theorems and $L^2$ estimates to solve the $\bar\partial$
equation for global holomorphic sections.  This additional copy of
the canonical line bundle has to be used to enable one to ignore
the effect of sufficiently small Lelong numbers.  Note that the
``two-tower argument'' for the deformational invariance of the
plurigenera also depends crucially on an additional copy of the
canonical line bundle to go down and up the two towers.

\bigbreak\noindent{\bf Part II. Application of Algebraic Geometry to
Analysis.}

\bigbreak We now consider the other direction which is the
application of algebraic geometry to subelliptic estimates.

\bigbreak\noindent{\sc (2.1) Regularity, Subelliptic Estimates,
and Kohn's Multiplier Ideals.}

\bigbreak Let us start out with the setting of the problem of
subelliptic estimates.

\bigbreak\noindent(2.1.1) {\it Setting of Regularity of the
$\bar\partial$-Problem}. Let $\Omega$ be a bounded domain in
${\mathbb C}^n$ with $C^\infty$ boundary $\partial\Omega$ which is
defined by a single $C^\infty$ real-valued function $r$ on some open
neighborhood of the topological closure $\bar\Omega$ of $\Omega$ in
${\mathbb C}^n$ so that $\Omega=\left\{r<0\right\}$ and
$\partial\Omega=\left\{r=0\right\}$ and the differential $dr$ of $r$
is nowhere zero on $\partial\Omega$.  We assume that
$\partial\Omega$ is weakly pseudoconvex in the sense that the
Hermitian form $\sqrt{-1}\partial\bar\partial r$ on
$T^{1,0}_{\partial\Omega}=T^{1,0}_{{\mathbb
C}^n}\cap\left(T_{\partial\Omega}\otimes{\mathbb C}\right)$ is
semipositive, where $T_{\partial\Omega}$ is the real tangent space
of the manifold $\partial\Omega$ and $T^{1,0}_{{\mathbb C}^n}$ is
the ${\mathbb C}$-vector space of all tangent vectors of type
$(1,0)$ of ${\mathbb C}^n$.

\medbreak\noindent(2.1.2) {\it Global Regularity Problem for
(0,1)-Forms.} The problem of global regularity of $\bar\partial$ for
$(0,1)$-forms asks whether, given a $C^\infty$ $\bar\partial$-closed
$(0,1)$-form $f$ on $\Omega$ which is $L^2$ with respect to the
Euclidean metric of ${\mathbb C}^n$, the solution $u$ of the
equation $\bar\partial u=f$ on $\Omega$ which is orthogonal to all
$L^2$ holomorphic functions on $\Omega$ with respect to the
Euclidean metric of ${\mathbb C}^n$ is $C^\infty$ up to the boundary
of $\Omega$.  The discussion here applies, with appropriate
adaptation, also to $(0,q)$-form for general $1\leq q\leq n-1$.  For
notational simplicity we will confine ourselves here only to the
case of $(0,1)$-forms.

\medbreak\noindent(2.1.3) {\it Subelliptic Estimates.}  The global
regularity problem is a consequence of the subelliptic estimate
which is defined as follows.  The {\it subelliptic estimate}, with
subellipticity order $\varepsilon>0$, holds at a boundary point $P$
of $\Omega$ if, for some open neighborhood $U_P$ of $P$ in ${\mathbb
C}^n$ and some positive number $C$,
$$
\left\|\left|\varphi\right|\right\|^2_\varepsilon\leq
C\left(\left\|\bar\partial\varphi\right\|^2+
\left\|\bar\partial^*\varphi\right\|^2+\left\|\varphi\right\|^2\right)
$$
for every test $(0,1)$-form $\varphi$ supported in
$U_P\cap\bar\Omega$ which is in the domain of $\bar\partial$ and in
the domain of $\bar\partial^*$, where $\left\|\cdot\right\|$ is the
usual $L^2$ norm on $\Omega$ and
$\left\|\left|\varphi\right|\right\|^2_\varepsilon$ is the Sobolev
$L^2$ norm on $\Omega$ for derivatives up to order $\varepsilon$ in
the direction tangential to $\partial\Omega$ ({\it i.e.,} directions
annihilated by $dr$).

\medbreak\noindent(2.1.4) {\it Kohn's Multipliers and their Ideals.}
To quantitatively measure the failure of subelliptic estimates Kohn
introduced multipliers [Ko79].  For a point $P$ in the boundary of
$\Omega$, a {\it Kohn multiplier} is a $C^\infty$ function germ $F$
on ${\mathbb C}^n$ at $P$ such that, for some open neighborhood
$U_P$ of $P$ in ${\mathbb C}^n$ and some positive numbers $C$ and
$\varepsilon$ (which may depend on $F$),
$$
\left\|\left|F\varphi\right|\right\|^2_\varepsilon\leq
C\left(\left\|\bar\partial\varphi\right\|^2+
\left\|\bar\partial^*\varphi\right\|^2+\left\|\varphi\right\|^2\right)
$$
for every test $(0,1)$-form $\varphi$ supported in
$U_P\cap\bar\Omega$ which is in the domain of $\bar\partial$ and in
the domain of $\bar\partial^*$.

\medbreak The set of all Kohn multipliers forms an ideal known as
the {\it Kohn multiplier ideal} which we denote by $I_P$.  By
definition it is clear that the subelliptic estimate holds if and
only if $I_P$ is the unit ideal (which means that the function
identically $1$ belongs to the ideal), in which case the order
$\varepsilon$ of subellipticity can be chosen to be the positive
number $\varepsilon$ for the identically $1$ function.

\medbreak Kohn [Ko79] introduced the following procedures of
generating multipliers.

\smallbreak\noindent(i)  $r\in I_P$.

\smallbreak\noindent(ii)  The coefficient of $\partial
r\wedge\bar\partial r\wedge\left(\partial\bar\partial
r\right)^{n-1}$ belongs to $I_P$.

\smallbreak\noindent(iii) $I_P$ equals its ${\mathbb R}$-radical
$\sqrt[{\mathbb R}]{I_P}$ in the sense that if $g\in I_P$ and
$\left|f\right|^m\leq\left|g\right|$ for some $m\geq 1$, then $f\in
I_P$.

\smallbreak\noindent(iv) If $f_1,\cdots,f_{n-1}$ belong to $I_P$ and
$1\leq j\leq n-1$, then the coefficients of
$$
\partial f_1\wedge\cdots\wedge\partial f_j\wedge
\partial r\wedge\bar\partial
r\wedge\left(\partial\bar\partial r\right)^{n-1-j}$$ belong to
$I_P$.

\bigbreak\noindent{\sc (2.2) Finite Type and Subellipticity.}

\bigbreak The problem is how to use geometric properties of the
boundary of $\Omega$ to conclude that the above procedures of Kohn
generate the unit ideal.

\medbreak\noindent(2.2.1) {\it Finite Type.} One natural geometric
property is the property of finite type which is defined as follows
[DA82, DK99].  The {\it type} $m$ at a point $P$ of the boundary of
$\Omega$ is the supremum of the normalized touching order
$$\frac{{\rm ord}_0\left(r\circ\varphi\right)}{{\rm ord}_0\varphi},$$
to $\partial\Omega$, of all local holomorphic curves
$\varphi:\Delta\to{\mathbb C}^n$ with $\varphi(0)=P$, where $\Delta$
is the open unit $1$-disk and ${\rm ord}_0$ is the vanishing order
at the origin $0$. A point $P$ of the boundary of $\Omega$ is said
to be of {\it finite type} if the type $m$ at $P$ is finite.  The
whole boundary $\partial\Omega$ is of finite type if the supremum
over the types of all points in the boundary of $\Omega$ is finite.

\medbreak\noindent(2.2.2) {\it Results of Kohn and Catlin.} Kohn's
original goal of developing his theory of multiplier ideals is to
show that his procedures of generating multipliers will yield the
unit ideal if $\partial\Omega$ is of finite type.

\medbreak For the case where $r$ is real-analytic, Kohn showed
[Ko79], by using a result of Diederich-Fornaess [DF78], that if his
multiplier ideal is not the unit ideal, it would lead to a
contradiction that $\partial\Omega$ contains some local holomorphic
curves.  His method of argument by contradiction does not yield an
effective $\varepsilon$ as an explicit function of the finite type
$m$ and the dimension $n$ of $\Omega$.

\medbreak Later Catlin [Ca83, Ca84, Ca87], using another approach of
multitypes and approximate boundary systems, showed that subelliptic
estimates hold for smooth bounded weakly pseudoconvex domains of
finite type.

\medbreak The problem remains whether Kohn's procedures are enough,
or some other procedures need to be found, to effectively generate
the unit ideal in the case of finite type.

\medbreak Catlin's use of multitype suggests the possibility of
additional procedures different from Kohn's which is geared to the
approach of multitype and which may be related to the ``Weierstrass
systems'' for ideals in the sense of Hironaka [AHV75, Hi77] or
Grauert [Ga73, Gr72].

\bigbreak\noindent(2.2.3) {\it Multiplier Modules.}  The listing of
Kohn's procedures of generating multipliers can be better organized
with the introduction of {\it multiplier modules}. For a point $P$
in the boundary of $\Omega$, a {\it multiplier-form} is a $C^\infty$
germ $\theta$ of $(1,0)$-form on ${\mathbb C}^n$ at $P$ such that,
for some open neighborhood $U_P$ of $P$ in ${\mathbb C}^n$ and some
positive numbers $C$ and $\varepsilon$ (which may depend on
$\theta$),
$$
\left\|\left|\theta\cdot\varphi\right|\right\|^2_\varepsilon\leq
C\left(\left\|\bar\partial\varphi\right\|^2+
\left\|\bar\partial^*\varphi\right\|^2+\left\|\varphi\right\|^2\right)
$$
for every test $(0,1)$-form $\varphi$ supported in
$U_P\cap\bar\Omega$ which is in the domain of $\bar\partial$ and in
the domain of $\bar\partial^*$, where $\theta\cdot\varphi$ is the
inner product between the $(1,0)$-form $\theta$ and the $(0,1)$-form
$\varphi$ with respect to the Euclidean metric of ${\mathbb C}^n$.
The {\it mutliplier module}, denoted by $A_P$, is the set of all
multiplier-forms at $P$. We break up Kohn's list of procedures into
the following three groups.

\medbreak\noindent(A) {\it Initial Membership}.

\smallbreak\noindent(i) $r\in I_P$.

\smallbreak\noindent(ii) $\partial\bar\partial_j r$ belongs to $A_P$
for every $1\leq j\leq n-1$ if $\partial r=\partial z_n$ at $P$ for
some local coordinate system $\left(z_1,\cdots,z_n\right)$, where
$\partial_j$ means $\frac{\partial}{\partial z_j}$.

\medbreak\noindent(B) {\it Generation of New Members}.

\smallbreak\noindent(i)  If $f\in I_P$, then $\partial f\in A_P$.

\smallbreak\noindent(ii) If $\theta_1,\cdots,\theta_{n-1}\in A_P$,
then the coefficient of
$$
\theta_1\wedge\cdots\wedge\theta_{n-1}\wedge\partial r$$ is in
$I_P$.

\medbreak\noindent(C) {\it Real Radical Property}.

\smallbreak\noindent If $g\in I_P$ and
$\left|f\right|^m\leq\left|g\right|$, then $f\in I_P$.

\bigbreak\noindent{\sc (2.3) Algebraic Formulation for Special
Domains.}

\bigbreak A {\it special domain} $\Omega$ in ${\mathbb C}^{n+1}$
(with coordinates $z_1,\cdots,z_n, w$) means a bounded domain
defined by
$$
{\rm
Re\,}w+\sum_{j=1}^N\left|h_j\left(z_1,\cdots,z_n\right)\right|^2<0,
$$
where $h_j\left(z_1,\cdots,z_n\right)$ is a holomorphic function
defined on some open neighborhood of the closure $\bar\Omega$ which
depends only on the first $n$ variables
$z=\left(z_1,\cdots,z_n\right)$.

\medbreak\noindent(2.3.1) {\it Finite Type Condition of Special
Domain.} The condition of finite type at a boundary point $P$ of
$\Omega$ can be formulated in terms of the ideal generated by
$h_1,\cdots,h_N$. To describe the formulation, we can assume without
loss of generality that $P$ is the origin of ${\mathbb C}^{n+1}$.
Let ${\mathcal O}_{{\mathbb C}^n,0}$ be the ring of all holomorphic
function germs of ${\mathbb C}^n$ at the origin of ${\mathbb C}^n$
and ${\mathfrak m}_{{\mathbb C}^n,0}$ be the maximum ideal of
${\mathcal O}_{{\mathbb C}^n,0}$.  Finite type at the origin means
that
$$
\left({\mathfrak m}_{{\mathbb
C}^n,0}\right)^p\subset\sum_{j=1}^N{\mathcal O}_{{\mathbb C}^n,0}h_j
$$
for some positive integer $p$.  The number $p$ is related to the
type $t$ in the following way.  The inequality
$$
\left|z\right|^q\leq C\sum_{j=1}^N\left|h_j(z)\right|
$$
holds for some positive constant $C$ on some open neighborhood of
the origin in ${\mathbb C}^n$ when $q\leq t\leq 2q$.  By Skoda's
theorem (1.8.1), $p$ can be chosen so that $p\leq q\leq(n+2)p$.

\medbreak\noindent(2.3.2) {\it Inductively Defined Ideals and
Functions.} We introduce the following inductively defined ideals
and positive-valued functions.

\medbreak For $\nu\in{\mathbb N}$ we inductively define ideals
$J_\nu$ and $\tilde J_\nu$ of ${\mathcal O}_{{\mathbb C}^n,0}$ and
positive-valued functions $\gamma_\nu$ on $J_\nu$ and
positive-valued functions $\tilde\gamma_\nu$ on $\tilde J_\nu$ as
follows.

\medbreak The ideal $J_1$ is generated by all elements $f$ such that
$f$ is defined by
$$
dg_1\wedge\cdots\wedge dg_n=f dz_1\wedge\cdots\wedge dz_n,
$$
where each of $g_1,\cdots,g_n$ is a ${\mathbb C}$-linear combination
of $h_1,\cdots,h_N$.

\medbreak The value $\gamma_1(f)$ is equal always to $\frac{1}{8}$
for $f\in J_1$. The ideal $\tilde J_1$ is the radical of $J_1$. For
$f\in\tilde J_1$ the value $\tilde\gamma_1(f)$ of the function
$\tilde\gamma_1$ at $f$ is equal to $\frac{1}{8m}$, where $m$ is the
smallest positive integer $m$ such that $f^m\in J_1$.

\medbreak For $\nu\in{\mathbb N}$ and $\nu\geq 2$ the ideal $J_\nu$
is generated by all elements $f$ such that $f$ is either an element
of $\tilde J_{\nu-1}$ or an element defined by an expression of the
form
$$
dg_1\wedge\cdots\wedge dg_n=f dz_1\wedge\cdots\wedge
dz_n,\leqno{(2.3.2.1)}
$$
where for some $0\leq k\leq n$ each of $g_1,\cdots,g_k$ is a
${\mathbb C}$-linear combination of $h_1,\cdots,h_N$ and each of
$g_{k+1},\cdots,g_n$ is an element of $\tilde J_{\nu-1}$.

\medbreak For $f\in J_\nu$ we assign the largest possible positive
number $a(f)$ which can be obtained in one of the following ways. If
$f$ is in $\tilde J_{\nu-1}$, the number $a(f)$ is the same as
$\tilde\gamma_{\nu-1}(f)$.  If $f$ is given by $(2.3.2.1)$,
\begin{itemize}\item[(i)]the number $a(f)$ is the minimum of
$$\frac{1}{8},\,\frac{1}{2}\,\tilde\gamma_{\nu-1}\left(g_{k+1}\right),\cdots,\frac{1}{2}\,\tilde\gamma_{\nu-1}\left(g_n\right)
$$ when $1\leq k<n\,$; and\item[(ii)] the number
$a(f)$ is the minimum of
$$\frac{1}{2}\,\tilde\gamma_{\nu-1}\left(g_1\right),\cdots,\frac{1}{2}\,\tilde\gamma_{\nu-1}\left(g_n\right)
$$ when $k=0\,$; and\item[(iii)] the number
$a(f)$ is equal to $\frac{1}{8}$ when $k=n$.\end{itemize}

\medbreak For $f\in\tilde J_\nu$ the value $\tilde\gamma_\nu(f)$ of
the function $\tilde\gamma_\nu$ at $f$ is equal to
$\frac{1}{m\gamma_\nu\left(f^m\right)}$, where $m$ is the smallest
positive integer $m$ such that $f^m\in J_\nu$.

\bigbreak The following statement is the algebraic formulation of
using Kohn's procedures to effectively generate the unit ideal in
the case of special domains of finite type.

\bigbreak\noindent(2.3.3) {\it Statement.} There exists a positive
number $\varepsilon$ which depends only on $n$ and $p$ by an
explicit expression such that for some $\nu\in{\mathbb N}$ there
exists some $f\in\tilde J_\nu$ with $f$ nonzero at the origin and
$\tilde\gamma_\nu(f)\geq\varepsilon$.

\medbreak From the point of view of local complex-analytic geometry,
the above statement (2.3.3) can be regarded as generalizing in a
very elaborate manner, to the case of ideals of holomorphic
functions of $n$ complex variables, the trivial statement that the
vanishing order of the differential of a holomorphic function germ
of a single complex variable at a point is equal to its vanishing
order minus $1$.  For this very elaborate generalization, the
process of taking differential is replaced by the process of taking
the Jacobian determinant of $n$ holomorphic function germs.

\bigbreak\noindent{\sc (2.4) Interpretation as Frobenius Theorem
over Artinian Subschemes.}

\bigbreak The problem of using Kohn's procedures to effectively
generate the unit ideal in the case of general type can be
interpreted as a Frobenius theorem over Artinian subschemes.

\medbreak The usual Frobenius theorem for ${\mathbb R}^m$ is the
following.  Let $U$ be an open subset of ${\mathbb R}^m$ open
subset.  Consider a smooth distribution $x\mapsto V_x\subset
T_{{\mathbb R}^m}={\mathbb R}^m$ of of $k$-dimensional subspaces of
$T_{{\mathbb R}^m}$.

\medbreak The Frobenius theorem can be formulation in terms of Lie
brackets or in terms of differential forms. The Lie bracket
formulation states that the distribution $V_x$ is integrable ({\it
i.e.} each $V_x$ is the tangent space of a submanifold of dimension
$k$ in a $(m-k)$-parameter family of such $k$-folds) if and only if
$V_x$ is closed under taking the Lie brackets of any two of its
elements ({\it i.e.,} the Lie bracket $\left[V_x,V_x\right]$ is
contained in $V_x$ for all $x\in U$).

\medbreak The differential form formulation state that the
distribution $V_x$ is integrable if and only if
$d\omega_j=\sum_{\ell=1}^{m-k}\omega_\ell\wedge\eta_\ell$ for some
$1$-forms $\eta_1,\cdots,\eta_{m-k}$, where
$\omega_1,\cdots,\omega_{m-k}$ are smooth $1$-forms defining $V_x$
({\it i.e.,} $V_x$ is the intersection of the kernels of
$\omega_1,\cdots,\omega_{m-k}$).

\medbreak By an Artinian subscheme we mean an unreduced subspace
supported at a single point.  In other words, it is a multiple
point.  For example, the ringed space $\left(0,\,{\mathcal
O}_{{\mathbb C}^n}\left/{\mathcal I}\right.\right)$, with
$$\left({\mathfrak m}_{{\mathbb C}^n,0}\right)^N\subset{\mathcal
I}$$ for some integer $N\geq 1$, is an Artinian subscheme.

\medbreak Now instead of considering integrability of the
distribution $V_x$ over the open subset $U$ of ${\mathbb R}^m$ we
will consider integrability over some multiple point.

\medbreak For our setting of bounded weakly pseudoconvex domain
$\Omega$ with smooth boundary $M=\partial\Omega$, the distribution
we consider is the space of all real tangent vectors of $M$ which
are the real parts of elements of $T_M^{(1,0)}$, where $T_M^{(1,0)}$
is the space of all complex-valued tangent vectors of $M$ of type
$(1,0)$. The Lie bracket formulation of the integrability of this
distribution over an open subset $M^\prime$ of $M$ is the same as
$M^\prime$ being Levi-flat.

\medbreak Integrability of this distribution over an Artinian
subscheme supported at a point $P$ of $M$ means the existence of a
local holomorphic curve touching $M$ at $P$ to an order
corresponding to the Artinian subscheme.

\medbreak The differential form formulation of the usual Frobenius
theorem now corresponds to the generation of new multipliers by
wedge product and exterior differentiation in Kohn's procedures.

\bigbreak\noindent{\sc (2.5) Sums of Squares and Kohn's
Counter-Example.}

\bigbreak In the case of complex dimension two one way in which Kohn
obtained subelliptic estiamtes for weakly pseudoconvex domains of
finite type [Ko72] is to use H\"ormander's subelliptic estimates for
sums of squares of real-valued vector fields whose iterated Lie
brackets span the entire tangent space [H\"o67].  The operator used
is $\bar\partial\bar\partial^*+\bar\partial^*\bar\partial$.  When
one decouples the operator from the domain and asks for
subellipticity from iterated Lie bracket conditions, one has to
consider the case of complex-valued vector fields in the setting of
sums of squares. Kohn recently produced the following
counter-example [Ko04].  Let $M$ be the boundary of the domain ${\rm
Re\,}z_2+\left|z_1\right|^2<0$ in ${\mathbb C}^2$ (with coordinates
$z_1, z_2$), which is biholomorphic to the complex $2$-ball, so that
$M$ is given by $ {\rm Re\,}z_2=-\left|z_1\right|^2$.  Let $x={\rm
Re\,}z_1$, $y={\rm Im\,}z_1$, $z=x+\sqrt{-1}y$, and $t={\rm
Im\,}z_2$.  Let
$$\displaylines{L=\frac{\partial}{\partial
z_1}-2\overline{z_1}\frac{\partial}{\partial
z_2}=\frac{\partial}{\partial z}+\sqrt{-1}\bar
z\frac{\partial}{\partial t},\cr\bar L=\frac{\partial}{\partial
\overline{z_1}}-2 z_1 \frac{\partial}{\partial
\overline{z_2}}=\frac{\partial}{\partial \bar z}-\sqrt{-1}
z\frac{\partial}{\partial t},\cr}$$ and
$$X_{1k}={\overline{z_1}\,}^k L,\quad X_2=\bar L, \quad
E_k=X_{1k}^*X_{1k}+X_2^*X_2.
$$
The commutators of $X_{1k}$, $X_2$ of order $\leq k+1$ span the
complexified tangent space of $M$ at $0$.  Kohn proved that for
$k>0$ subellipticity does not hold for $E_k$, yet hypoellipticity
holds for $E_k$ with a loss of $k$ derivatives in the supremum norm
and a loss of $k-1$ derivatives in the Sobolev norm for $k>1$.

\medbreak\noindent(2.5.1) {\it Explanation of Failure of
Subellipticity.} The reason for the failure of subellipticity can
be illustrated by the following domain version corresponding to
the boundary version in Kohn's counter-example.  Consider the ball
$\Omega$ defined by $\left|z_1-1\right|^2+\left|z_2\right|^2-1<0$
in ${\mathbb C}^2$ and the vector field
$$
L=\overline{z_2}\frac{\partial}{\partial
z_1}+\left(1-\overline{z_1}\right)\frac{\partial}{\partial z_2}
$$
which spans the space $T^{1,0}_{\partial\Omega}$ of complex-valued
tangent vectors of $\partial\Omega$.  Let
$$
X_1=\overline{z_1}L,\quad X_2=\bar L.
$$
Then the commutators of $X_1$, $X_2$ of order $\leq 2$ span the
complexified tangent space of $\partial\Omega$ at $0$.  Let
$E=X_1^*X_1+X_2^*X_2$.  Unlike $\bar\partial\bar\partial^*+
\bar\partial^*\bar\partial$ which corresponds to the situation with
$X_1$ replaced by $L$, the operator $E$ does not have
subellipticity.

\medbreak Take any $\varepsilon>0$. To see the failure of
subellipticity for $E$ with subellipticity order $\varepsilon>0$,
take $p>2$ and $C>0$ so that, if $g$ is a function on $\Omega$ with
$L^2$ norm $\leq 1$ and if its Sobolev $L^2$ norm on $\Omega$ for
derivative up to order $\varepsilon$ is also $\leq 1$, then the
$L^p$ norm of $g$ on $\Omega$ is $\leq C$.  For $0<\eta<1$ and
$0<\alpha<1$, take a branch of
$\varphi_\eta=\frac{1}{\left(z_1+\eta\right)^\alpha}$ on $\Omega$.
The $L^2$ norm of $\varphi_\eta$ on $\Omega$ is bounded uniformly in
$0<\eta<1$.  Choose $\alpha$ so close to $1$ that the $L^p$ norm of
$\varphi_\eta$ on $\Omega$ is not bounded uniformly in $0<\eta<1$.

\medbreak Because of the factor $\overline{z_1}$ in $X_1$, the $L^2$
norm of $X_1\varphi_\eta$ on $\Omega$ is bounded uniformly in
$0<\eta<1$. Since $X_2\varphi_\eta$ is identically zero, if the
subelliptic estimate holds for $E$ with subellipticity order
$\varepsilon>0$, the Sobolev norm of $\varphi_\eta$ on $\Omega$ for
derivative up to $\varepsilon$ would be bounded uniformly in
$0<\eta<1$ and, as a consequence, the $L^p$ norm of $\varphi_\eta$
on $\Omega$ would be bounded uniformly for $0<\eta<1$, which
contradicts the choice of $\alpha$.

\bigbreak\noindent(2.5.2) {\it Sums of Squares of Matrix-Valued
Vector Fields.}  The subellipticity of
$\bar\partial\bar\partial^*+\bar\partial^*\bar\partial$ for
$(0,1)$-forms holds on weakly pseudoconvex domains of finite type.
The operator $\bar\partial\bar\partial^*+\bar\partial^*\bar\partial$
for $(0,1)$-forms can be written as a sum of squares of
matrix-valued vector fields.  On the other hand, even for the
strongly pseudoconvex case of the complex $2$-ball, when the
operator is decoupled from the domain and different from
$\bar\partial\bar\partial^*+\bar\partial^*\bar\partial$, Kohn's
counter-example shows that, unlike the case of {\it real-valued}
vector fields, in general subellipticity fails for a sum of squares
of {\it complex-valued} vector fields with the iterated Lie bracket
condition. An important natural problem is to understand what
additional conditions are involved in the case of
$\bar\partial\bar\partial^*+\bar\partial^*\bar\partial$ for
$(0,1)$-forms on weakly pseudoconvex domains of finite type which
would give subellipticity for sums of squares of matrix-valued
vector fields in general.

\bigbreak\noindent{\it References.}

\medbreak\noindent [AHV75] J. M. Aroca, H. Hironaka, and J. L.
Vicente, The theory of the maximal contact. {\it Memorias de
Matemática del Instituto ``Jorge Juan''}, No. 29. Instituto ``Jorge
Juan'' de Matem\'aticas, Consejo Superior de Investigaciones
Cientificas, Madrid, 1975.

\medbreak\noindent[AS95] U. Angehrn and Y.-T. Siu, Effective
freeness and point separation for adjoint bundles. {\it Invent.
Math.} 122 (1995), 291--308.

\medbreak\noindent[BT76] E. Bedford and B. A. Taylor, The
Dirichlet problem for a complex Monge-Amp\'ere equation. {\it
Invent. Math.} 37 (1976), 1--44.

\medbreak\noindent[Ca83] D. Catlin, Necessary conditions for the
subellipticity of the $\bar\partial$-Neumann problem, {\it Ann. of
Math.} 117 (1983), 147-171.

\medbreak\noindent[Ca84] D. Catlin, Boundary invariants of
pseudoconvex domains. {\it Ann. of Math.} 120 (1984), 529--586.

\medbreak\noindent[Ca87] D. Catlin, Subelliptic estimates for the
$\bar\partial$-Neumann problem on pseudoconvex domains. {\it Ann.
of Math.} 126 (1987), 131-191.

\medbreak\noindent[DA82] J.~D'Angelo, Real hypersurfaces, orders
of contact, and applications. {\it Ann. of Math.} 115 (1982),
615--637.

\medbreak\noindent[DK99] J.~D'Angelo and J.~J.~Kohn, Subelliptic
estimates and finite type. Several complex variables (Berkeley,
CA, 1995--1996), 199--232, {\it Math. Sci. Res. Inst. Publ.} 37,
Cambridge Univ. Press, Cambridge, 1999.

\medbreak\noindent[De93] J.-P. Demailly, A numerical criterion for
very ample line bundles.  {\it J. Diff. Geom.} 37 (1993),
323--374.

\medbreak\noindent[De96a] J.-P. Demailly, $L\sp 2$ vanishing
theorems for positive line bundles and adjunction theory. In {\it
Transcendental methods in algebraic geometry} (Cetraro, 1994),
1--97, {\it Lecture Notes in Math.} 1646, Springer, Berlin, 1996.

\medbreak\noindent[De96b] J.-P. Demailly, Effective bounds for
very ample line bundles. {\it Invent. Math.} 124 (1996), no. 1-3,
243--261.

\medbreak\noindent[DF78] K. Diederich and J. E. Fornaess,
Pseudoconvex domains with real-analytic boundary. {\it Ann. Math.}
107 (1978), 371--384.

\medbreak\noindent[EL93] L. Ein and R. Lazarsfeld, Global generation
of pluricanonical and adjoint linear series on smooth projective
threefolds. {\it J. Amer. Math. Soc.} 6 (1993), 875--903.

\medbreak\noindent[ELN94] L. Ein, R. Lazarsfeld, and M. Nakamaye,
Zero-estimates, intersection theory, and a theorem of Demailly.
{\it Higher-dimensional complex varieties} (Trento, 1994),
183--207,

\medbreak\noindent[FdB96] G. Fernández del Busto, Matsusaka-type
theorem on surfaces. {\it J. Alg. Geom.} 5 (1996), 513--520.

\medbreak\noindent[Fu87] T. Fujita, On polarized manifolds whose
adjoint bundles are not semipositive, Proceedings of the 1985 Sendai
Conference on Algebraic Geometry, {\it Advanced Studies in Pure
Mathematics} 10, 167-178 (1987).

\medbreak\noindent[Fu93] T. Fujita, Remarks on Ein-Lazarsfeld
criterion of spannedness of adjoint bundles of polarized
threefolds, preprint, 1993.

\medbreak\noindent[Ga73] A. Galligo, `A propos du th\'eor\`eme de
pr\'eparation de Weierstrass. {\it Fonctions de plusieurs variables
complexes} (S\'em. Fran\c{c}ois Norguet, octobre 1970--d\'ecembre
1973, \`a la m\'emoire d'Andr\'e Martineau), pp. 543--579. {\it
Lecture Notes in Math.} Vol. 409, Springer, Berlin, 1974.

\medbreak\noindent[Gr72] H. Grauert, \"Uber die Deformation
isolierter Singularit\"aten analytischer Mengen. {\it Invent. Math.}
15 (1972), 171--198.

\medbreak\noindent[Hei02] G. Heier, Effective freeness of adjoint
line bundles. {\it Doc. Math.} 7 (2002), 31--42

\medbreak\noindent[Hel97] S. Helmke, On Fujita's conjecture. {\it
Duke Math. J.} 88 (1997), 201--216.

\medbreak\noindent[Hel99] S. Helmke, On global generation of
adjoint linear systems. {\it Math. Ann.} 313 (1999), 635--652.

\medbreak\noindent[Hi64] H. Hironaka, Resolution of singularities of
an algebraic variety over a field of characteristic zero. I, II.
{\it Ann. of Math.} 79 (1964), 109--203; {\it ibid.} 79 (1964),
205--326.

\medbreak\noindent[Hi77] H. Hironaka, Bimeromorphic smoothing of a
complex-analytic space. {\it Acta Math. Vietnam} 2 (1977), 103--168.

\medbreak\noindent[H\"o67] L. H\"ormander, Hypoelliptic second order
differential equations. {\it Acta Math.} 119 (1967), 147--171.

\medbreak\noindent[Ka82] Y. Kawamata, A generalization of the
Kodaira-Ramanujam's vanishing theorem, {\it Math. Ann.} 261
(1982), 43--46.

\medbreak\noindent[Ka85] Y. Kawamata, Pluricanonical systems on
minimal algebraic varieties. {\it Invent. Math.} 79 (1985),
567--588.

\medbreak\noindent[Ka97] Y. Kawamata, On Fujita's freeness
conjecture for $3$-folds and $4$-folds. {\it Math. Ann.} 308
(1997), 491--505.

\medbreak\noindent[Ko72] J. J. Kohn, Boundary behavior of
$\bar\partial$ on weakly pseudo-convex manifolds of dimension two.
Collection of articles dedicated to S. S. Chern and D. C. Spencer on
their sixtieth birthdays. {\it J. Diff. Geom.} 6 (1972), 523--542.

\medbreak\noindent[Ko79] J. J. Kohn, Subellipticity of the
$\bar\partial$-Neumann problem on pseudo-convex domains: sufficient
conditions. {\it Acta Math.} 142 (1979), 79-122.

\medbreak\noindent[Ko04] J. J. Kohn, Hypoellipticity and loss of
derivatives. {\it Ann. of Math.} (to appear), Preprint 2004.

\medbreak\noindent[Le83] M. Levine, Pluri-canonical divisors on
K\"ahler manifolds, {\it Invent. math.} 74 (1983), 293--903.

\medbreak\noindent[Le85] M. Levine, Pluri-canonical divisors on
K\"ahler manifolds, II. {\it Duke Math. J.} 52 (1985), no. 1,
61--65.

\medbreak\noindent[Mo88] S. Mori, Flip theorem and the existence of
minimal models for $3$-folds. {\it J. Amer. Math. Soc.} 1 (1988),
117--253.

\medbreak\noindent[Nad89] A. Nadel, Multiplier ideal sheaves and
existence of K\"ahler-Einstein metrics of positive scalar curvature,
{\it Proc. Nat. Acad. Sci. U.S.A.} 86, 7299-7300 (1989) and {\it
Ann. of Math.} 132, 549-596 (1990).

\medbreak\noindent[Nak86] N. Nakayama, Invariance of the plurigenera
of algebraic varieties under minimal model conjectures. {\it
Topology} 25 (1986), 237--251.

\medbreak\noindent[OT87] T. Ohsawa and K. Takegoshi, On the
extension of $L^2$ holomorphic functions, {\it Math. Zeitschr.}
195 (1987), 197--204.

\medbreak\noindent[Re88] I. Reider, Vector bundles of rank 2  and
linear systems on algebraic surfaces, {\it Ann. of Math.} 127,
309-316 (1988).

\medbreak\noindent[Si74] Y.-T. Siu, Analyticity of sets associated
to Lelong numbers and the extension of closed positive currents.
{\it Invent. Math.} 27 (1974), 53--156.

\medbreak\noindent[Si93] Y.-T. Siu, An effective Matsusaka big
theorem. {\it Ann. Inst. Fourier} (Grenoble) 43 (1993), 1387--1405.

\medbreak\noindent[Si94]  Y.-T. Siu, Very ampleness criterion of
double adjoints of ample line bundles. In {\it Modern methods in
complex analysis} (Princeton, NJ, 1992), 291--318, Ann. of Math.
Stud. 137, Princeton Univ. Press, Princeton, NJ, 1995

\medbreak\noindent[Si96a]  Y.-T. Siu, The Fujita conjecture and
the extension theorem of Ohsawa-Takegoshi. In {\it Geometric
Complex Analysis} ed. Junjiro Noguchi {\it et al}, World
Scientific: Singapore, New Jersey, London, Hong Kong 1996, pp.
577--592.

\medbreak\noindent[Si96b] Y.-T. Siu, Effective very ampleness.
{\it Invent. Math.} 124 (1996), 563--571.

\medbreak\noindent[Si98] Y.-T. Siu, Invariance of plurigenera,
{\it Invent. Math.} 134 (1998), 661-673.

\medbreak\noindent[Si01] Y.-T. Siu, Very ampleness part of Fujita's
conjecture and multiplier ideal sheaves of Kohn and Nadel.  In: {\it
Complex Analysis and Geometry}, ed. J. D. McNeal, Volume 9 of the
Ohio State University Mathematical Research Institute Publications,
Verlag Walter de Gruyter, 2001, pp.171-191.

\medbreak\noindent[Si02a] Y.-T. Siu, Extension of twisted
pluricanonical sections with plurisubharmonic weight and invariance
of semipositively twisted plurigenera for manifolds not necessarily
of general type. In: {\it Complex Geometry (Collection of papers
dedicated to Hans Grauert)}, ed. I. Bauer {\it et al}, Springer,
Berlin, 2002, pp. 223--277.

\medbreak\noindent[Si02b] Y.-T. Siu, A New Bound for the Effective
Matsusaka Big Theorem, {\it Houston J. Math.} 28 (issue for 90th
birthday of S.S. Chern) (2002), 389-409.

\medbreak\noindent[Si03] Y.-T. Siu, Invariance of Plurigenera and
Torsion-Freeness of Direct Image Sheaves of Pluricanonial Bundles,
In: {\it Finite or Infinite Dimensional Complex Analysis and
Applications} (Proceedings of the 9th International Conference on
Finite or Infinite Dimensional Complex Analysis and Applications,
Hanoi, 2001), edited by Le Hung Son, W. Tutschke, C.C. Yang, Kluwer
Academic Publishers 2003, pp.45-84.

\medbreak\noindent[Sk72] H. Skoda, Application des techniques $L^2$
\`a la th\'eorie des ideaux d'un alg\`ebre de fonctions holomorphes
avec poids, {\it Ann. Sci. Ec. Norm. Sup.} 5 (1972), 548-580.

\medbreak\noindent[Vi82] E. Viehweg, Vanishing theorems. {\it J.
Reine Angew. Math.} 335 (1982), 1--8.

\medbreak\noindent[Ya78] S.-T. Yau, On the Ricci curvature of a
compact Kähler manifold and the complex Monge-Amp\`ere equation. I.
{\it Comm. Pure Appl. Math.} 31 (1978), 339--411.

\bigbreak\noindent{\it Author's Address:} Department of Mathematics,
Harvard University, Cambridge, MA 02138

\smallbreak\noindent{\it e-mail:} siu@math.harvard.edu

\end{document}